\documentclass[oneside, 12pt]{amsart}
\usepackage{amssymb,mathrsfs,amstext,enumerate,comment}
\usepackage[noadjust]{cite}
\usepackage[plainpages=false,colorlinks,hyperindex,pdfpagemode=UseNone,bookmarksopen,linkcolor=blue,citecolor=blue,urlcolor=blue]{hyperref}

\newtheorem{theorem}{Theorem}[section]
\newtheorem{proposition}[theorem]{Proposition}
\newtheorem{corollary}[theorem]{Corollary}
\newtheorem{lemma}[theorem]{Lemma}

\theoremstyle{definition}

\newtheorem{example}[theorem]{Example}

\theoremstyle{remark}

\newtheorem{problem}[theorem]{Problem}

\numberwithin{equation}{section}

\def\Span{{\rm span}\,}

\def\tr{{\rm tr}\,}

\def\dim{{\rm dim}\,}

\def\IC{{\mathbb C}}
\def\IF{{\mathbb F}}
\def\IR{{\mathbb R}}

\def\bM{{\mathbf M}}
\def\bU{{\bf U}}

\newcommand{\TT}{{\mathrm{t}}}

\parskip=\baselineskip

\begin{document}
\openup 1\jot

\title[Norm-parallelism preservers]{Linear maps of \texorpdfstring{$C^*$}{C*}-algebras preserving norm parallelism}

\author[Li, Tsai, Wang and Wong]{Chi-Kwong Li, Ming-Cheng Tsai, Ya-Shu Wang \and Ngai-Ching Wong}

\address[Li]{Department of Mathematics, The College of William \& Mary, Williamsburg, VA 23187, USA.}
\email{ckli@math.wm.edu}

\address[Tsai]{General Education Center, National Taipei University of Technology, Taipei 10608, Taiwan.}
\email{mctsai2@mail.ntut.edu.tw}

\address[Wang]{Department of Applied Mathematics, National Chung Hsing University, Taichung 40227, Taiwan.}
\email{yashu@nchu.edu.tw}

\address[Wong]{Department of Applied Mathematics, National Sun Yat-sen University, Kaohsiung 80424, Taiwan; Department of Healthcare Administration and Medical Information, Kaohsiung Medical University, Kaohsiung 80708, Taiwan.}
\email{wong@math.nsysu.edu.tw}

\thanks{Corresponding author: Ngai-Ching Wong, wong@math.nsysu.edu.tw}

\date{October 14, 2026}

\begin{abstract}
Two elements $a,b$ in a normed space are parallel if
$\|a+\gamma b\|=\|a\|+\|b\|$ for some $|\gamma|=1$;
they form a triangle equality attaining (TEA) pair if
$\|a+b\|=\|a\|+\|b\|$.
For finite-dimensional complex $C^*$-algebras other than
$\mathbf M_2(\mathbb C)$ and $\mathbb C\oplus_\infty\mathbb C$, bijective linear
maps on these algebras preserving parallel or TEA pairs
are scalar multiples of isometries.
For unital abelian algebras with spectra containing at least
three points, preservation in both directions gives the same
conclusion.
We obtain corresponding results for type I factors, simple
type I $C^*$-algebras, and unitarily invariant operator spaces
under additional continuity and finite-rank hypotheses.
All results for $C(X)$ and $B(H,K)$ hold over both real
and complex fields.
\end{abstract}



\maketitle

\section{Introduction}\label{s:Intro}

The geometry of the unit sphere of a $C^*$-algebra reflects its algebraic
structure. Kadison's classical theorem \cite{Kad51} provides a fundamental example
of this connection: every surjective complex-linear isometry between
unital $C^*$-algebras is a Jordan $^*$-isomorphism followed by multiplication
by a unitary element. This raises a natural question: how much of the geometry of
the unit sphere must a linear bijection preserve in order to determine
the same algebraic structure? In this paper, we consider a condition
expressing the local flatness of the unit sphere, namely preservation of
norm parallelism.

Two elements $a,b$ in a normed space $E$ over $\mathbb F=\mathbb R$ or
$\mathbb C$ are said to be \emph{norm parallel}, or simply
\emph{parallel}, if
\[
\|a+\lambda b\|=\|a\|+\|b\|
\qquad\text{for some $\lambda\in\mathbb F$ with $|\lambda|=1$.}
\]
In this case, we write $a\parallel b$. If this equality holds specifically with $\lambda=1$, then
$(a,b)$ is called a \emph{triangle equality attaining} (TEA) pair.
For nonzero $a,b$, the Hahn--Banach theorem shows that $(a,b)$ is a TEA
pair if and only if there exists a norm-one functional $\varphi\in E^*$
such that
\[
\varphi(a)=\|a\|\qquad\text{and}\qquad\varphi(b)=\|b\|.
\]
Equivalently, the line segment joining $a/\|a\|$ and $b/\|b\|$ lies
entirely on the unit sphere $S_E$. Accordingly, $a\parallel b$ if and
only if, for some unimodular scalar $\lambda$, the segment joining
$a/\|a\|$ and $\lambda b/\|b\|$ lies in $S_E$.
Thus parallelism records flat portions of the unit sphere after a scalar
phase adjustment. A parallel pair of linearly independent vectors gives
a nondegenerate segment of this kind.

The significance of this relation depends on the geometry of the space.
In a strictly convex normed space, two nonzero vectors are parallel
exactly when they are linearly dependent; they form a TEA pair exactly
when one is a positive scalar multiple of the other. Consequently, every
linear bijection between strictly convex normed spaces preserves both
relations. For spaces with non-strictly convex norms, however, parallelism
contains additional geometric information. Our purpose is to show that,
for several major classes of $C^*$-algebras and operator spaces, preserving this
information forces a linear bijection to preserve the entire norm up to
a constant. In the complex $C^*$-algebra setting, the underlying Jordan
$^*$-structure is then completely determined up to multiplication by a unitary (in the multiplier algebra in the nonunital case).

Related preserver problems have been studied for function spaces and
matrix spaces. For $L^1$ spaces, bijective linear preservers of parallel or TEA pairs are characterized by boundedness and disjointness preservation
\cite[Theorems 4.2 and 4.3]{LTWW-lp}; by Kan's theorem \cite[Theorem 4.1]{Kan}, such maps admit classical Lamperti-type representations via
regular set homomorphisms and weight functions.
Related results on disjointness, triple products, and norm preservers
between rectangular matrix spaces appear in \cite{LTWW}.
For rectangular matrix spaces equipped with the operator norm,
\cite{LTWW-matrix} gives a complete classification. When both matrix
dimensions are at least two, and the space is not $\bM_2(\mathbb F)$,
every bijective linear parallelism preserver is a scalar multiple of a
standard isometry. Further studies of norm parallelism and triangle
equality can be found in
\cite{BCM19,Ket,Maligranda08,Manna27,NT02,Popovici12,Wojcik,ZM15,ZM16}.

Our finite-dimensional result provides a direct connection between local
flatness and the structure of a $C^*$-algebra. If $\mathcal A$ is a
finite-dimensional complex $C^*$-algebra other than
$\mathbb C\oplus_\infty\mathbb C$ and $\bM_2(\mathbb C)$, then
the following conditions on a bijective linear map
$T:\mathcal A\to\mathcal A$ are equivalent:
\begin{enumerate}[\quad(1)]
 \item $T$ preserves parallel pairs;
 \item $T$ preserves TEA pairs;
 \item $T$ is a scalar multiple of an isometry.
\end{enumerate}
In particular, preservation in one direction suffices; no preservation
assumption on $T^{-1}$ is needed
(Theorem~\ref{thm:finite-Calg}). Equivalently, such a map has the form
\[
T(a)=\gamma uJ(a)\qquad(a\in\mathcal A),
\]
where $\gamma>0$, $u$ is a unitary in $\mathcal A$, and $J$ is a Jordan
$^*$-automorphism of $\mathcal A$.
The result is obtained from a classification on finite
$\ell_\infty$-sums of rectangular matrix spaces with both block dimensions
at least two, together with scalar summands
(Theorem~\ref{thm:Sum-R}). The direct-sum setting introduces a further
rigidity: preservation of parallelism identifies the matrix summands and
forces a common scale on all components. In particular, the presence of
a $\mathbb C\oplus_\infty\mathbb C$ or $\bM_2(\mathbb C)$ summand does not itself produce an exception
when the algebra has additional nonzero summands.

For abelian unital $C^*$-algebras, we obtain an infinite-dimensional
counterpart without assuming continuity. Let $X,Y$ be compact Hausdorff
spaces with at least three points. A linear bijection $T:C(X)\to C(Y)$
preserves parallelism in both directions if and only if
\[
Tf(y)=u(y)f(\sigma(y))\qquad(f\in C(X),\ y\in Y),
\]
where $\sigma:Y\to X$ is a homeomorphism and $u\in C(Y)$ has constant
nonzero modulus (Theorem~\ref{thm:C(X)}). The same characterization holds
for preservation of TEA pairs in both directions. Thus the relation
describing local flatness recovers both the underlying compact space
and the norm up to scale. The proof reconstructs points from families
of functions attaining their norms there and uses the rigidity of
linear peak-set preservers (Theorem~\ref{thm:C(X)-Peak}).

We also study parallelism preservers on unitarily invariant standard
subspaces of $B(H,K)$, where $H,K$ are real or complex Hilbert spaces.
Here a subspace is standard if it contains every finite-rank operator.
For surjective linear maps preserving parallelism in both directions,
Theorem~\ref{main-2} gives scalar multiples of the standard unitary or
transpose isometries under either of the following additional hypotheses:
the map is weak operator continuous and sends finite-rank operators to
finite-rank operators; or its domain consists of compact operators and
it maps the finite-rank operators onto the finite-rank operators.
These hypotheses yield corresponding results for type I factors and
compact-operator algebras
(Corollaries~\ref{cor:typeIfactors} and~\ref{cor:K(H)}).

The matrix and operator arguments examine the linear spans of the
parallelism sets
\[
\mathcal P(A)=\{B:B\parallel A\}.
\]
Their codimensions, together with their weak operator closures in the
infinite-dimensional setting, reveal the maximal-singular-value
structure of $A$. This allows preservation of the geometric relation
to be converted into preservation of disjoint rank-one operators,
from which the standard representations follow.

Every linear TEA-pair preserver preserves parallel pairs, but the
converse need not hold. The distinction is already visible in the
exceptional complex spaces $\mathbb C\oplus_\infty\mathbb C$ and
$\bM_2(\mathbb C)$: their bijective parallelism preservers need
not be scalar multiples of isometries, whereas their bijective
complex-linear TEA preservers are of this form. We recall the relevant
classifications in Lemma~\ref{lem:row-col} and
Propositions~\ref{prop-parallel-C} and~\ref{prop-parallel-R}, including
the different behavior over the real field.
For finite-dimensional real normed spaces, a linear bijection preserves parallel pairs if and only if it preserves TEA pairs \cite[Corollary 2.19]{Manna27}.

These results lead to the following general question concerning the extent to
which local flatness determines the structure of a $C^*$-algebra.

\begin{problem}\label{prob:1}
For which complex $C^*$-algebras $\mathcal A,\mathcal B$ must every
linear bijection $T:\mathcal A\to\mathcal B$ preserving parallelism
in both directions be a scalar multiple of an isometry?
Does the same conclusion hold for arbitrary complex $C^*$-algebras
if TEA pairs are preserved in both directions?
\end{problem}

Throughout this paper, function spaces $C(X)$, Hilbert spaces $H,K$, and rectangular operator spaces in $B(H,K)$ may
be over either $\mathbb{R}$ or $\mathbb{C}$, with linearity understood over the indicated field.
On the other hand, all $C^*$-algebras are assumed to be complex. All operator spaces are equipped with the operator
norm.

The paper is organized as follows. Section~\ref{s:cont} treats continuous-function spaces,
Section~\ref{s:fc} finite-dimensional $C^*$-algebras and rectangular
matrix sums, and Section~\ref{s:B(H,K)} standard operator spaces.

\section{Linear parallelism preservers of continuous functions}\label{s:cont}

In this section, we study the linear parallelism preservers between abelian unital $C^*$-algebras,
i.e., the algebras of continuous functions on a compact Hausdorff space.
All results and their proofs hold for both the real and the complex cases.
For simplicity, we assume the underlying field $\IF$ is the complex field $\mathbb{C}$.
If $\IF=\IR$, then we can simply skip those statements involving non-real numbers, and the arguments
will remain valid.

Let $X$ be a compact Hausdorff space.
Let $f\in C(X)$ be a continuous function on $X$.
Denote by
$$
\operatorname{Pk}(f) = \{x\in X: |f(x)| = \|f\|\}
$$
the \emph{peak set} of $f$, which consists of points at which $f$ assumes its supnorm.
Clearly, $\operatorname{Pk}(\mathbf{1})=\operatorname{Pk}(\mathbf{0})=X$, where
$\mathbf{1}$ and $\mathbf{0}$ are the constant one and zero functions on $X$, respectively.
Indeed, $\operatorname{Pk}(f)=X$ exactly when the absolute value $|f|$ of
$f$ is a constant function, that is,
$|f|=\alpha \mathbf{1}$ for some scalar $\alpha$.

Two continuous functions $f,g\in C(X)$ are parallel (denote by $f\parallel g$) exactly when they have a common
peak point $x$, that is, $|f(x)|=\|f\|$ and $|g(x)|=\|g\|$. In other words,
\begin{align}\label{eq:f||g}
  f \parallel g \quad\text{if and only if}\quad \operatorname{Pk}(f)\cap\operatorname{Pk}(g)\neq \emptyset\qquad
  \text{for every $f,g\in C(X)$.}
\end{align}
Recall that $f,g\in C(X)$ are said to be \emph{disjoint} if $fg=0$.
For nonzero disjoint functions $f,g$, we have
$\operatorname{Pk}(f)\cap \operatorname{Pk}(g)=\emptyset$, and thus disjoint nonzero functions cannot be parallel.

The following lemma describes
 the structure of linear parallelism preservers of finite dimensional abelian $C^*$-algebras
 $\ell_\infty^n:= \underbrace{\IF\oplus_\infty \cdots \oplus_\infty \IF}_n$. Recall that a generalized permutation matrix is a square matrix in which
there is exactly one nonzero entry in each row and each column, which has modulus one.

\begin{lemma}[{\cite[Theorem 2.6]{LTWW-lp}}]\label{lem:row-col}
   Let  $T: {\ell_\infty^n}\to {\ell_\infty^n}$ be a linear map.
 \begin{enumerate}
 \item [{\rm (a)}]
 $T$ preserves parallel pairs if and only if there is $\gamma \ge 0$ and a generalized permutation matrix $Q$
 such that one of the following forms holds:

 \begin{enumerate}
  \item[{\rm (a.1)}] $T$ has the form $\mathbf{x} \mapsto \gamma Q\mathbf{x}$,
  \item[{\rm (a.2)}] $n = 2$ and $T$ has the form $\mathbf{x} \mapsto \gamma CQ\mathbf{x}$, where
  $C =\begin{pmatrix} 1 & \beta \\ \overline{\beta} & 1\end{pmatrix}$
for some scalar $\beta$ with $|\beta|<1$.
  \item[{\rm (a.3)}] $T=\mathbf{v}{\mathbf{u}}^\TT$ for some nonzero column vectors ${\mathbf{u}}, {\mathbf{v}}\in \IF^n$.
\end{enumerate}
\item[{\rm (b)}]  $T$  preserves
TEA pairs if and only if
one of the following holds.

\begin{enumerate}
\item[{\rm (b.1)}] $T$ has the form in {\rm (a.1)}.

\item[{\rm (b.2)}] $\IF^n = \IR^2$ and $T$ has the form in {\rm (a.2)}.

\item[{\rm (b.3)}] $\IF^n = \IR^2$ and $T$ has the form in {\rm (a.3)}, where
$\mathbf{u}=(u_1, u_2)^{\TT}$ with $|u_1| = |u_2|$.
\end{enumerate}
 \end{enumerate}
\end{lemma}

We also need the following result to establish Theorem \ref{thm:C(X)}. It is also
of independent interest.

\begin{theorem}\label{thm:C(X)-Peak}
  Let $X$ be a compact Hausdorff space containing at least three points.
Let $T: C(X)\to C(X)$ be a linear map preserving peak sets, namely,
$$
 \operatorname{Pk}(Tf)=\operatorname{Pk}(f) \quad\text{for every $f\in C(X)$.}
$$
Then $u=T\mathbf{1}\in C(X)$ has constant nonzero modulus, i.e., $|u|=\lambda\mathbf{1}$
for some scalar $\lambda>0$, and
$$
Tf=uf \quad\text{for every $f\in C(X)$.}
$$
\end{theorem}
\begin{proof}
We divide the proof into several steps.

  \textbf{Step 1.}
  Since $\operatorname{Pk}(T\mathbf{1})=\operatorname{Pk}(\mathbf{1})=X$, there is a non-negative scalar
  $\lambda$ such that $|T\mathbf{1}|=\lambda \mathbf{1}$. If $\lambda=0$ then
  $$
 \operatorname{Pk}(f)= \operatorname{Pk}(Tf)=\operatorname{Pk}(Tf + T\mathbf{1})=\operatorname{Pk}(T(f+\mathbf{1}))=\operatorname{Pk}(f+\mathbf{1}) \quad\text{for every $f\in C(X)$.}
 $$
 This is impossible, for some negative functions $f$, say,
  and thus $\lambda\neq 0$. Replacing $T$ with $(T\mathbf{1})^{-1} T$, we can
 assume that $u:=T\mathbf{1}=\mathbf{1}$.

 Suppose $Tf=\mathbf{0}$ for some $f\in C(X)$. Then $\operatorname{Pk}(f)=\operatorname{Pk}(Tf)=X$ implies that $|f|=\alpha\mathbf{1}$ for some scalar
 $\alpha\geq 0$. Without loss of generality,
  we assume that $f(x_0)=\alpha$ for some point $x_0\in X$.
 Since $T$ is unital,
 $$
 \operatorname{Pk}(f+ \alpha\mathbf{1}) =\operatorname{Pk}(T(f+ \alpha\mathbf{1}))
 = \operatorname{Pk}(\alpha\mathbf{1})=X.
 $$
 This gives
 $$
 |f(x) + \alpha| = f(x_0) + \alpha = 2\alpha= 2\|f\| \quad\text{for all $x\in X$.}
 $$
It follows $f=\alpha\mathbf{1}$. In this case,
  $\mathbf{0}=Tf=\alpha T\mathbf{1}
=\alpha\mathbf{1}$ implies $\alpha=0$, and thus, $f=\mathbf{0}$.
  This shows that $T$ is injective.

  From now on, we assume that $T$ is a unital and injective linear map preserving peak sets.

  \textbf{Step 2.}
Let $x_1, x_2, x_3$ be three distinct points in $X$. Let $f_1, f_2, f_3\in C(X)$ be nonnegative
continuous functions with norm one
such that
$$
f_1f_2=f_2f_3=f_3f_1=\mathbf{0}\quad\text{and}\quad f_1(x_1)=f_2(x_2)=f_3(x_3)=1.
$$
Let $d_{jk}=Tf_j(x_k)$ for $j,k=1,2,3$.
We claim that the matrix $D=(d_{jk})$ is the $3\times 3$ identity matrix.

We first note that $Tf_j$ peaks at $x_j$ but not other $x_k$ as $f_j$ does.
Hence,
$$
|d_{jj}|> |d_{jk}|\quad\text{ whenever $j\neq k$.}
$$
Since $f_1 + f_2 + t\mathbf{1}$ and $f_1 + f_2 + \mathrm{i}t\mathbf{1}$ peak at $x_1, x_2$, we see that
$Tf_1 + Tf_2 + t\mathbf{1}$ and $Tf_1 + Tf_2 + \mathrm{i}t\mathbf{1}$ also peak at $x_1, x_2$ for any positive $t$.
It follows
$$
|d_{11} + d_{21} + t| = |d_{12} + d_{22} + t| \quad\text{and}\quad |d_{11} + d_{21} + \mathrm{i}t| = |d_{12} + d_{22} + \mathrm{i}t|,
$$
or
$$
|d_{11} + d_{21}|^2 + 2t\operatorname{Re}(d_{11} + d_{21}) + t^2
= |d_{12} + d_{22}|^2 + 2t\operatorname{Re}(d_{12} + d_{22}) + t^2
$$
and
$$
|d_{11} + d_{21}|^2 + 2t\operatorname{Im}(d_{11} + d_{21}) + t^2
= |d_{12} + d_{22}|^2 + 2t\operatorname{Im}(d_{12} + d_{22}) + t^2,
$$
for any $t>0$.
This forces
$$
d_{11} + d_{21}= d_{12} + d_{22}.
$$
Exchanging the roles of $f_1, f_2, f_3$, we get
$$
d_{11} + d_{31}=d_{13}+d_{33} \quad\text{and}\quad d_{22}+d_{32}=d_{23}+d_{33}.
$$

Consider $f_1 + f_2 + f_3 + t\mathbf{1}$ and $f_1 + f_2 + f_3 + \mathrm{i}t\mathbf{1}$ for $t>0$.
The same argument gives
$$
d_{11} + d_{21} + d_{31} = d_{12} + d_{22} + d_{32}=d_{13}+d_{23}+d_{33}.
$$
This forces
$$
d_{31}= d_{32}, \quad d_{21}=d_{23}\quad\text{and}\quad  d_{12}=d_{13}.
$$

Since both $f_1 + f_2 + s f_3$ and $f_1 - f_2 + s f_3$ peak at $x_1, x_2$ for any scalar $s$ with $|s|\leq 1$,
 we have
$$
|d_{11}\pm d_{21}+ sd_{31}| = |  d_{12}\pm d_{22}+ sd_{32}|,
$$
and thus
$$
|d_{11}\pm d_{21}|^2 + 2\operatorname{Re}\, [sd_{31}\overline{(d_{11}\pm d_{21})}] + |sd_{31}|^2   =
|d_{12}\pm d_{22}|^2 + 2\operatorname{Re}\, [sd_{32}\overline{(d_{12}\pm d_{22})}] + |sd_{32}|^2.
$$
If $d_{31}=d_{32}\neq 0$, because $s$ is arbitrary but $|s|\leq 1$, we have
$$
d_{11}\pm d_{21} =   d_{12}\pm d_{22}.
$$
This forces $d_{11}=d_{12}$, which conflicts with the fact $|d_{11}|>|d_{12}|$.
Hence, $d_{31}=d_{32}=0$. Similarly, we have $d_{12}=d_{13}=0$ and $d_{21}=d_{23}=0$.
It also follows $d_{11}=d_{22}=d_{33}$.

Suppose $d_{11} = \alpha + i\beta$ for $\alpha, \beta\in \mathbb{R}$.
Since $f_1 - f_2 - i\beta\mathbf{1}$ peaks at $x_1, x_2$, we have
$$
|d_{11} - d_{21} - i \beta| = |d_{12} - d_{22}  -i\beta|,
$$
or
$$
|\alpha| = |\alpha + 2i\beta|.
$$
Hence $\beta=0$.
Since $f_1 - \frac{1}{2}\mathbf{1}$ peaks at $x_1$ and $x_2$,
 and so does $Tf_1 - \frac{1}{2}\mathbf{1}$, we have
 $$
 \left|d_{11} - \frac{1}{2}\right| = \left|d_{12} - \frac{1}{2}\right|,
\quad\text{or}\quad
\left|\alpha - \frac{1}{2}\right| = \frac{1}{2}.
$$
The last equality gives $\alpha\in\{0,1\}$. Since $T$ is injective and
$x_1$ is a peak point of $Tf_1$, $|\alpha|=\|Tf_1\|>0$.
Therefore, $d_{11}=\alpha=1$.
This establishes the claim that $D$ is the $3\times 3$ identity matrix.

\textbf{Step 3.}
We first prove that $T$ is bounded. Step 2 has the following consequence:
if $g\geq0$ vanishes on some neighborhood of two distinct points $p,q$, then
\[
 \|Tg\|=\|g\|,\qquad Tg(p)=Tg(q)=0.
\]
Indeed, for $g\ne0$, choose nonnegative norm-one functions supported in
disjoint neighborhoods of $p,q$ on which $g$ vanishes, and apply Step 2
to these functions and $g/\|g\|$.

Fix three distinct points $p_1,p_2,p_3\in X$. The open sets
$U_i=X\setminus\{p_j:j\ne i\}$ cover $X$. Choose a nonnegative continuous
partition of unity $\rho_1,\rho_2,\rho_3$ subordinate to this cover,
with $\operatorname{supp}\rho_i\subseteq U_i$.
For every $f\geq0$, each $\rho_i f$ vanishes on neighborhoods of the other
two points. Consequently,
\[
 \|Tf\|\leq\sum_{i=1}^3\|T(\rho_i f)\|
 =\sum_{i=1}^3\|\rho_i f\|\leq3\|f\|.
\]
Decomposing a real function into its positive and negative parts gives
$\|Tf\|\leq6\|f\|$. In the complex case, decomposing the real and imaginary
parts gives $\|Tf\|\leq12\|f\|$. Thus $T$ is bounded over either field.

\textbf{Step 4.}
We show that $T$ is local. Fix $x\in X$ and suppose first that $g\geq0$
vanishes on a neighborhood of $x$. Choose distinct $p,q\in X\setminus\{x\}$
and a nonnegative partition of unity $\eta_p+\eta_q=1$ subordinate to
$X\setminus\{p\}$ and $X\setminus\{q\}$. Each of $\eta_p g,\eta_q g$
vanishes near $x$ and one further point. The consequence of Step 2 used
above therefore gives $Tg(x)=0$. Taking positive and negative parts, and
then real and imaginary parts if necessary, gives the same conclusion
for every $g$ vanishing near $x$.

If $f(x)=0$, for every $\epsilon>0$ choose an open neighborhood $U$ of $x$
on which $|f|<\epsilon$, and a continuous function $0\leq\chi\leq1$
which vanishes near $x$ and equals one on $X\setminus U$.
Then $\|f-\chi f\|\leq\epsilon$ and $T(\chi f)(x)=0$.
Boundedness implies $Tf(x)=0$. Applying this to $f-f(x)\mathbf1$ and using
$T\mathbf1=\mathbf1$ yields $Tf(x)=f(x)$ for every $f$ and every $x$.
Undoing the normalization proves the theorem.
\end{proof}

While the case when $X, Y$ consist of at most two points are done in Lemma \ref{lem:row-col}, we have

\begin{theorem}\label{thm:C(X)}
Let $X,Y$ be compact Hausdorff spaces containing at least three points.
Let $T: C(X)\to C(Y)$ be a bijective linear map. The following are equivalent.
\begin{enumerate}[(a)]
  \item Both $T, T^{-1}$ preserve parallel pairs.
  \item Both $T, T^{-1}$ preserve TEA pairs.
  \item There is a homeomorphism $\sigma: Y\to X$ and a function $u\in C(Y)$ with constant modulus
such that
$$
Tf(y) = u(y)f(\sigma(y)) \quad\text{for any $f\in C(X)$ and $y\in Y$.}
$$
\end{enumerate}
\end{theorem}
\begin{proof}
It is plain that (c) $\implies$ (b) $\implies$ (a). Suppose (a) holds, and we verify (c).
  For any $x\in X$, let
  $$
  {\operatorname{Pk}^{-1}_x} = \{ f\in C(X): |f(x)|=\|f\|\}
  $$
to be the family of continuous functions $f\in C(X)$ peaking at $x$.
Consequently, the peak set of any $f\in C(X)$ is the compact set
$$
\operatorname{Pk}(f)=\{x\in X: f\in {\operatorname{Pk}^{-1}_x}\}.
$$
The same notions are defined for any point $y\in Y$ and for
any $g\in C(Y)$.
Let
$$
{\mathcal{Q}_x} = \bigcap \{\operatorname{Pk}(Tf): f\in {\operatorname{Pk}^{-1}_x}\}.
$$

We claim that ${\mathcal{Q}_x}$ is nonempty.
Suppose, otherwise, by compactness of $Y$,
there are finitely many nonzero
$f_1,\ldots, f_n\in {\operatorname{Pk}^{-1}_x}$ such that $\operatorname{Pk}(Tf_1) \cap \cdots\cap \operatorname{Pk}(Tf_n) =\emptyset$.
Let $V_j = Y\setminus \operatorname{Pk}(Tf_j)$ for $j=1,\ldots, n$. Then $Y = V_1 \cup \ldots \cup V_n$ is an open covering of $Y$ such that
$Tf_j$ does not peak in $V_j$ for $j=1,\ldots, n$.
We can assume that $f_1(x)=\cdots=f_n(x)=1$.
Let $f= \min_{1\leq j \leq n} |f_j|\in {\operatorname{Pk}^{-1}_x}$.
Let the nonzero continuous function $Tf$ peaks at some $y\in Y$.
We can assume that $y\in V_1$.
Let $g_1\in C(X)$ such that $Tg_1\in {\operatorname{Pk}^{-1}_y}$ and $Tg_1$ vanishes outside $V_1$.
Since $Tf_1, Tg_1$ do not have a common peak point, they are not parallel with each other, and so are not $f_1, g_1$.
Since $\operatorname{Pk}(f)\subseteq \operatorname{Pk}(f_1)$, we see that $f, g_1$, and thus $Tf, Tg_1$ are not parallel with each other.
This contradicts the fact that $|Tf(y)|=\|Tf\|$ and $|Tg_1(y)|=\|Tg_1\|$.
Hence, ${\mathcal{Q}_x}$ is nonempty.

Interchanging the roles of $T$ and $T^{-1}$, we see that
$$
{\mathcal{Q}_y}=\bigcap \{\operatorname{Pk}(g): Tg\in {\operatorname{Pk}^{-1}_y}\} \neq \emptyset \quad\text{for any $y\in Y$.}
$$
Let $y\in {\mathcal{Q}_x}$ and $z\in {\mathcal{Q}_y}$. In other words,
for any $f\in C(X)$, if $f$ peaks at $x$ then $Tf$ peaks at $y$, and $f=T^{-1}(Tf)$ peaks at $z$.
It follows $z=x$, and thus ${\mathcal{Q}_y}=\{x\}$ is a singleton. Similarly, ${\mathcal{Q}_x}=\{y\}$ is also a singleton.
This gives us a bijective map $\sigma: Y\to X$ such that
\begin{quote}
$f\in C(X)$ peaks at $x=\sigma(y)$ if and only if $Tf$ peaks at $y$.
\end{quote}

We will see that both $\sigma$ and $\sigma^{-1}$ are continuous. For example, if $y_\lambda\to y$ in $Y$
but $\sigma(y_\lambda)\not\to \sigma(y)$ in the compact space
$X$, by replacing $\{y_\lambda\}_\lambda$ with a subnet we can choose
 an $f\in C(X)$ such that $f$ peaks at
all $\sigma(y_\lambda)$ but vanishes at $\sigma(y)$. Then $Tf$ peaks at all $y_\lambda$ but not at $y$. This contradicts the
continuity of $Tf$. Hence, $\sigma$ is continuous. Similarly, $\sigma^{-1}$ is continuous, and thus
$\sigma$ is a homeomorphism, from $Y$ onto $X$.

Composing $T$ with the isomorphism from $C(Y)$ onto $C(X)$ defined by $g\mapsto g\circ \sigma^{-1}$, we can assume that $X=Y$ and
$\sigma(x)=x$ for all $x\in X$.
In this setting, $\operatorname{Pk}(Tf)=\operatorname{Pk}(f)$ for any $f\in C(X)$.
Then Theorem \ref{thm:C(X)-Peak} applies and gives the desired assertion.
\end{proof}

\section{Linear parallelism preservers of finite dimensional \texorpdfstring{$C^*$}{C*}-algebras}\label{s:fc}

In this section, we discuss linear parallel/TEA pair preservers between finite dimensional (complex)
$C^*$-algebras.
Recall that every finite dimensional $C^*$-algebra
is an $\ell_\infty$ direct sum of full complex matrix algebras,
i.e., $\bM_{n_1}\oplus_\infty \cdots \oplus_\infty \bM_{n_k}\oplus_\infty \underbrace{\IC \oplus_\infty  \cdots \oplus_\infty  \IC }_h$,
 where all $n_i\geq 2$.
Our result states:

\begin{theorem}\label{thm:finite-Calg}
  Let $\mathcal{A}$ be a finite dimensional $C^*$-algebra such that $\mathcal{A}\neq\bM_{2}(\IC)$
  and $\mathcal{A}\neq\IC\oplus_\infty \IC$.
  The following statements are equivalent for a bijective linear map $T: \mathcal{A}\to\mathcal{A}$.
  \begin{enumerate}[(a)]
    \item $T$ preserves parallel pairs.
    \item $T$ preserves TEA pairs.
    \item $T$ is a scalar multiple of an isometry.
   \end{enumerate}
\end{theorem}

In a broader sense, we consider such linear preservers between
a (real or complex linear) $\ell_\infty$-direct
sum of rectangular matrices
$$
{\bM} := \bM_{m_1,n_1}\oplus_\infty  \cdots \oplus_\infty  \bM_{m_k,n_k}
\oplus_\infty  \underbrace{\IF \oplus_\infty  \cdots \oplus_\infty  \IF }_h,
$$
where the underlying field $\IF=\IR$ or $\IC$.
The structure theory of linear parallel/TEA pair preservers of finite dimensional $C^*$-algebras will
follow as a special case of such preservers of rectangular matrices.

Note that
the norms of the row matrix space $\bM_{1, n}$ and the column matrix space $\bM_{m,1}$ both
reduce to the $\ell_2$-norm of $\IF^n$ and $\IF^m$,
which are strictly convex.
Consequently, two nonzero row or column matrices
  $A, B$ in $\bM_{1, n}$ or $\bM_{m,1}$ are parallel (resp.\ form a TEA pair)
  exactly when $A=\lambda B$ for some scalar $\lambda$ (resp.\ with $\lambda >0$).
  Thus \emph{any} linear map of row matrix spaces or
    column matrix spaces preserves parallel pairs as well as TEA pairs.

In view of the following
propositions, we see that linear parallel/TEA pair preservers
of (real or complex) rectangular matrices are merely
scalar multiples of isometries, except for the $2\times 2$ matrix case.

\begin{proposition}[{\cite[Theorem 1.1]{LTWW-matrix}}] \label{prop:matrices}
Suppose $m,n \ge 2$, $(m,n) \ne (2,2)$,
and $T:\bM_{m,n}(\IF)\rightarrow \bM_{m,n}(\IF)$ is a nonzero linear map.
 The following conditions are equivalent.
\begin{enumerate}[{\rm (a)}]
\item $T$ preserves TEA pairs.
\item $T$ preserves parallel pairs and the range space of $T$ has dimension larger than one.
\item There are $U \in \bU_m(\IF)$, $V \in \bU_n(\IF)$, and $\gamma > 0$ such that

\medskip\hspace{0pt}{
{\rm (1)} \ $T$ has the form $A \mapsto \gamma UAV$, \ \ or \ \
{\rm (2)} \ $m = n$ and $T$ has the form $A \mapsto \gamma UA^{\TT}V$.
}
\end{enumerate}

If the range space of a linear map $T: \bM_{m,n}(\IF) \rightarrow \bM_{m,n}(\IF)$ has dimension
at most $1$, or equivalently,
$T$ has the form

\hspace{20pt}{
{\rm (3)} \
$A \mapsto \tr(FA) Z$ \quad for some $F\in \bM_{n,m}(\IF)$ and $Z\in \bM_{m,n}(\IF)$,
}

\noindent
then clearly $T$ preserves parallel pairs (but not necessarily TEA pairs).
\end{proposition}

The $2\times 2$ case appears to be much more complicated, in which
such preservers need \emph{not} be scalar multiples of
  isometries as demonstrated in Propositions \ref{prop-parallel-C} and \ref{prop-parallel-R}.
Let
$$
I_2 = \begin{pmatrix} 1 & 0 \cr 0 & 1\cr\end{pmatrix}, \
\mathcal{C}_1 = \mathcal{R}_1 = \begin{pmatrix} 1 & 0 \cr 0 & -1\cr\end{pmatrix}, \
\mathcal{C}_2 = \mathcal{R}_2 = \begin{pmatrix} 0 & 1 \cr 1 & 0\cr\end{pmatrix}, \
\mathcal{C}_3 = \begin{pmatrix} 0 & -\mathrm{i} \cr \mathrm{i} & 0 \cr\end{pmatrix}, \
\mathcal{R}_3 = \begin{pmatrix} 0 & 1 \cr -1 & 0 \cr\end{pmatrix}.
$$
It is well-known that
\begin{itemize}
\item $\mathbf{C} = \{I_2/\sqrt 2, \mathcal{C}_1/\sqrt 2, \mathcal{C}_2/\sqrt 2, \mathcal{C}_3/\sqrt 2\}$
is an orthonormal basis of $\bM_2(\IC)$
 with respect to the complex inner product
$\langle A,B\rangle=\operatorname{trace}(B^*A)$, while
\item
 $\mathbf{R}=\{I_2/\sqrt 2, \mathcal{R}_1/\sqrt 2, \mathcal{R}_2/\sqrt 2, \mathcal{R}_3/\sqrt 2\}$
is an orthonormal basis of $\bM_2(\IR)$ with respect to the real inner product $\langle A, B\rangle = \operatorname{trace}(B^{\TT} A)$.
\end{itemize}

\begin{proposition}[{\cite[Theorem 1.2]{LTWW-matrix}}]\label{prop-parallel-C}
Let $T: \bM_2(\IC) \rightarrow \bM_2(\IC)$ be a complex linear map.
\begin{itemize}
\item[(a)]$T$ preserves TEA pairs if and only if there are $\gamma \geq 0$ and
$U, V \in \bU_2(\IC)$ such that $T$ has the form
$$
A \mapsto \gamma UAV \qquad \hbox{ or } \qquad  A\mapsto \gamma UA^{\TT}V.
$$
\item[(b)] $T$ preserves parallel pairs if and only if one of the following holds.
\begin{itemize}
\item[(b.1)] $T$ has the form $A \mapsto \tr(FA) Z$ for some $F, Z\in \bM_2(\IC)$.
\item[(b.2)] There are complex unitary matrices
 $U_1, U_2, V_1, V_2 \in \bU_2(\IC)$ such that the matrix representation of the linear transformation
  $\tilde{T}$ defined by $A\mapsto U_1 T(U_2AV_2)V_1$ with respect
to the basis $\mathbf{C}$ has the form
\begin{equation*}
\ \hskip .7 in [\tilde{T}]_{\mathbf{C}} =\begin{pmatrix}
d_0 & \gamma_1 i & \gamma_2 i & \gamma_3 i \\
0 & d_1 & 0 & 0 \cr
0 & 0 & d_2 & 0 \cr
0 & 0 & 0 & d_3\cr
\end{pmatrix}
\end{equation*}
for real scalars
$d_0, d_1, d_2, d_3, \gamma_1, \gamma_2, \gamma_3$.
\end{itemize}
\end{itemize}
\end{proposition}

\begin{proposition}[{\cite[Theorem 1.3]{LTWW-matrix}}]\label{prop-parallel-R}
Let $T: \bM_2(\IR) \rightarrow \bM_2(\IR)$ be a real linear map.
\begin{itemize}
\item[(a)] $T$ preserves TEA pairs if and only if $T$ has the form:
\begin{itemize}
\item[(a.1)] there are real orthogonal matrices $F_1, F_2 \in \bU_2(\IR)$ with $F_1^{\,\TT}F_2 = \begin{pmatrix} 0 & 1\\ -1 & 0\end{pmatrix}$
and $Z_1, Z_2\in \bM_2(\IR)$ such that $T$
has the form
$$X \mapsto \tr(F_1X)Z_1 + \tr(F_2 X)Z_2;$$
\item[(a.2)] there are real orthogonal matrices $U_1, U_2, V_1, V_2 \in \bU_2(\IR)$
such that the matrix representation of the linear transformation
  $\tilde{T}$ defined by $A\mapsto U_1 T(U_2AV_2)V_1$ with respect to the basis $\mathbf{R}$ has the diagonal form
\begin{align*}
[\tilde T]_\mathbf{R} = \begin{pmatrix} d_0 & & &\\ & d_1 & &\\ & & d_2 &\\ & & & d_3\end{pmatrix}\in \bM_4(\IR).
\end{align*}
\end{itemize}
\item[(b)] $T$ preserves parallel pairs if and only if $T$ has the form in (a.1), (a.2), or
$$
A\mapsto \tr(FA)Z \quad\text{for some $F, Z \in \bM_2(\IR)$.}
$$
\end{itemize}
\end{proposition}

We are ready to state the main result in this section, from which it follows
Theorem \ref{thm:finite-Calg} directly.
Indeed, the form \eqref{eq:Tdecom} completely describes the structure
of linear parallel/TEA pair preservers of finite dimensional $C^*$-algebras other than
$\bM_2(\IC)$ and $\IC\oplus_\infty \IC$.

\begin{theorem}\label{thm:Sum-R}
Let the underlying field be $\mathbb{F}=\mathbb{R}$ or $\mathbb{C}$.
Let ${\bM} = \bM_{m_1,n_1}\oplus_\infty  \cdots \oplus_\infty  \bM_{m_k,n_k}
\oplus_\infty  \underbrace{\IF \oplus_\infty  \cdots \oplus_\infty  \IF }_h$ be an
$\ell_\infty$ direct sum, in which $m_j, n_j\geq 2$,
  for $j=1,\ldots, k$.
  Let $T: {\bM}\to {\bM}$
 be a bijective linear map.
If $\bM \neq \IF\oplus_\infty\IF$ and $\bM\neq \bM_2(\IF)$ then
 the following statements are equivalent.
 \begin{enumerate}[{\rm (a)}]
 \item $T$ preserves parallel pairs.
 \item $T$ preserves TEA pairs.
  \item $T$ is a scalar multiple of an isometry.
\item There exist a positive
  scalar $\gamma=\|T\|$,
unitary matrices $U_i\in\bM_{m_i}$ and $V_i\in\bM_{n_i}$,
permutations $\sigma$ of $\{1,\ldots,k\}$ and $\tau$ of $\{1,\ldots,h\}$,
and unimodular scalars $u_1,\ldots,u_h$, such that
\begin{equation}\label{eq:Tdecom}
 T\left(\bigoplus_{j=1}^k A_j\oplus(a_1,\ldots,a_h)\right)
 =\gamma\left(\bigoplus_{i=1}^k U_i A_{\sigma(i)}^{\dag_i}V_i
 \oplus(u_1a_{\tau(1)},\ldots,u_ha_{\tau(h)})\right).
\end{equation}
Here $\dag_i$ is either the identity, in which case
$(m_{\sigma(i)},n_{\sigma(i)})=(m_i,n_i)$, or transposition, in which case
$(n_{\sigma(i)},m_{\sigma(i)})=(m_i,n_i)$.

\end{enumerate}
\end{theorem}

For the exceptional cases,
the structure of linear parallel/TEA pair preservers of $\bM=\IF\oplus_\infty \IF$ and $\bM=\bM_2(\IF)$ are described
in Lemma \ref{lem:row-col} and Propositions \ref{prop-parallel-C} and \ref{prop-parallel-R}, respectively.

We might assume $k\geq 1$, for else it will reduce
 to the $\ell_\infty^h\cong \underbrace{\IF \oplus_\infty  \cdots \oplus_\infty  \IF }_h$
  case in Lemma \ref{lem:row-col}.
Moreover, the case when $\bM=\bM_{2,n}$ with $n\geq 3$ is covered in Proposition \ref{prop:matrices}.
We can thus assume
$$
m_1+\cdots + m_k+h\geq 3.
$$
The implications (d) $\Longrightarrow$ (c) $\Longrightarrow$ (b) $\Longrightarrow$ (a)
are plain.

 We introduce the following notations to facilitate our proof for (a) $\Longrightarrow$ (d).
For $A \in \bM_{m,n}$, let
$$
\mathcal{P}_{\bM_{m,n}}(A) =\{B\in \bM_{m,n}: B \text{ is parallel with }  A\}.
$$
In a more general setting,
for
$$
A =A_1\oplus\cdots\oplus A_k\oplus (a_1) \oplus\cdots\oplus (a_h)
\in {\bM} = \bM_{m_1,n_1}\oplus_\infty \cdots \oplus_\infty \bM_{m_k,n_k}
\oplus_\infty \underbrace{\IF \oplus_\infty \cdots \oplus_\infty \IF }_h,
$$
let
\begin{align*}
\mathcal{P}_{\bM}(A) = \{B =\bigoplus_{j=1}^k B_j\oplus (b_1) \oplus\cdots\oplus (b_h)\in {\bM}: B \text{ is parallel with }  A\},
\end{align*}
and let
\begin{quote}
\quad $\Span\mathcal{P}_{\bM}(A)$ be the linear span of all elements in $\bM$ parallel with $A$.
\end{quote}

With $m,n\geq 2$,
every real or complex rectangular matrix
$A \in \bM_{m,n}$ has a singular value decomposition; namely, there are unitary
 matrices
$U \in \bM_m$ and $V \in \bM_n$ such that
$U^*AV = \sum_{j=1}^r s_j(A) E_{jj}$, where $r$ is the rank of $A$, and
$\|A\| =s_1(A) \ge \cdots \ge s_r(A)$ are the positive singular values of $A$ in descending order.
Here, depending on the underlying field $\IF$ being the complex $\IC$ or the real $\IR$,
 by a unitary matrix we mean either a complex unitary matrix or a
   real orthogonal matrix.
More precisely, if $u_1, \dots, u_m$ are the columns of $U$ and $v_1, \dots, v_n$ are
the columns of $V$, then
\begin{equation}\label{eq:SVD}
A = U(\sum_{j=1}^r s_j(A)E_{jj})V^* =
\sum_{j=1}^r s_j(A) u_j v_j^* = \sum_{j=1}^r s_j(A) u_j\otimes v_j
\end{equation}
is a positive linear sum of disjoint rank one partial isometries $u_j\otimes v_j:= u_jv_j^*$.
 Here, we say that two rectangular matrices $B,C\in \bM_{m,n}$ are \emph{disjoint}, or \emph{orthogonal}
  as defined in \cite{BFGMP08},
 if $B^*C= \mathbf{0}_n$ and $BC^*=\mathbf{0}_m$.

A rank one partial isometry $u\otimes v$ is said to be
a \emph{peaking partial isometry} of $A$ if $|\langle Av, u\rangle| = \|A\|$.
In particular, any $u_j\otimes v_j$ in \eqref{eq:SVD} with
$s_j(A)=\|A\|$ is a peaking partial isometry of $A$.
Two rectangular matrices
$A, B \in \bM_{m,n}$ are parallel (resp.\ form a TEA pair)
 if and only if there are unitary matrices $U\in \bM_m,V \in \bM_n$
so that
$U^*AV=\begin{bmatrix}a\end{bmatrix}\oplus A_1$ and $U^*BV=\begin{bmatrix}b\end{bmatrix}\oplus B_1$
with $\|A\|=|a|$ and
$\|B\|=|b|$ (resp.\ $\|A\|=a$ and $\|B\|=b$); see, e.g.,
\cite{Li0}. In other words, $A, B$ are parallel exactly when they
share a common peaking isometry $u\otimes v= uv^*$.

We may regard $A=A_1\oplus\cdots\oplus A_k\oplus (a_1) \oplus\cdots\oplus (a_h)$
 in ${\bM}$ as a diagonal block matrix
in $\bM_{m+h,n+h}$, where
$$
m=m_1+\cdots + m_k\quad\text{and}\quad n=n_1+\cdots + n_k.
$$
The singular value decompositions of the direct components $A_1$, \ldots, $A_k$, together with its
diagonal part in $\underbrace{\IF \oplus \cdots \oplus \IF }_h$,
provide us a singular value decomposition of $A$. In particular, the peaking partial isometries in this blockwise
singular value decomposition can be chosen to lie in a block $\bM_{m_j, n_j}$ for some $j=1,\ldots,k$, or appears in the diagonal part $\underbrace{\IF \oplus \cdots \oplus \IF }_h$.

Let
\begin{align*}
\operatorname{Pr}_p =\{A\in \bM: \ & \text{the singular values of $A$
listed in descending order satisfy} \\
&\hspace{3cm} s_1(A) = s_2(A) = \cdots =s_p(A) > s_{p+1}(A) \geq \cdots\}.
\end{align*}
We call matrices in $\operatorname{Pr}_p$ matrices with \emph{peak rank} $p$.

\begin{lemma} \label{lem:Pk1}
Let $\bM = \bigoplus_{j=1}^k \bM_{m_j,n_j}\oplus
\underbrace{\IF \oplus_\infty \cdots \oplus_\infty \IF }_h$ be an
$\ell_\infty$ direct sum, in which the underlying field
 $\IF=\IR$ or $\IC$, and $m_j, n_j \geq 2$ for $j=1,\ldots, k$.
Let $A =\bigoplus_{j=1}^k A_j\oplus (a_1)\oplus\cdots\oplus(a_h)\in {\bM}$.
\begin{itemize}
\item[{\rm (1)}]
If $A$ has peak rank greater than one, i.e., $s_1(A) = s_2(A)$, then $\Span \mathcal{P}_\bM(A)$ has codimension at most $1$ in ${\bM}$.
\item[{\rm (2)}]
Suppose $A$ has peak rank one, i.e., $s_1(A) > s_2(A)$. There are two exclusive cases.
\begin{enumerate}[{\rm (a)}]
\item There are unitary matrices $U_i \in \bM_{m_i}, V_i \in \bM_{n_i}$
such that $U_i^*A_iV_i = [a] \oplus F$ with $a=\|A\| =s_1(A) > \|F\|, \|A_j\|, |a_l|$
for some index $i$, all $j\neq i$ and $l=1,\ldots,h$. Then
\begin{align*}
\hspace{1cm}\mathcal{P}_\bM(A)&=  \Bigg\{\bigoplus_j G_j\oplus (g_1) \oplus\cdots\oplus (g_h)\in \bM: G_i = U_i ([g] \oplus G)V_i^*\ \text{with}\  G \in \bM_{m_i-1,n_i-1}, \\
& G_j\in \bM_{m_j,n_j}, g, g_l\in \IF,\ \text{such that}\ |g| \geq \|G\|, \|G_j\|,    |g_l| \ \text{for}\ j\neq i, l=1,\ldots, h \Bigg\},\ \text{and}\\
 \Span&\mathcal{P}_\bM(A)=  \Bigg\{\bigoplus_j B_j\oplus (b_1) \oplus\cdots\oplus (b_h)\in \bM: B_i = U_i ([b] \oplus F)V_i^*\\
&\text{with arbitrary}\
 F \in \bM_{m_i-1,n_i-1},  B_j\in \bM_{m_j,n_j},\ \text{and}\ b, b_1,\ldots, b_h\in \IF\Bigg\}.
\end{align*}
In this case,
$\Span \mathcal{P}_{\bM}(A)$ has codimension $m_i + n_i -2\geq 2$ in ${\bM}$.

\item There is $1\leq l \leq h$ such that
 $|a_l| = \|A\| > \|A_i\|, |a_j|$ for $i=1,\ldots, k$ and $j\neq l$. In this case,
$\Span \mathcal{P}_\bM (A)=\bM$.
\end{enumerate}

\item[{\rm (3)}]
If $\Span \mathcal{P}_{\bM}(A)$ has codimension greater than one, then $A$ has peak rank one
and the codimension of $\Span \mathcal{P}_{\bM}(A) = m_i + n_i-2$ in $\bM$ for some index
$i=1,\ldots, k$.
In this case, condition {\rm (2a)} holds.
\end{itemize}
\end{lemma}
\begin{proof}
(1) If $s_1(A_i) = s_2(A_i)=\|A\|$ for some component rectangular
matrix $A_i\in\bM_{m_i,n_i}$, we may assume that $\|A\|=\|A_i\|=1$ and
$A_i = I_2
\oplus Z$ for some $Z \in \bM_{m_i-2,n_i-2}$ with $\|Z\| \le 1$.
Then $\mathcal{P}_{\bM_{m_i,n_i}}(A_i)$ contains matrices of the form
$[b] \oplus X$ such that $X\in \bM_{m_i-1,n_i-1}$
with $\|X\|\leq |b|$.
It follows that $B=(b_{ij})\in \Span \mathcal{P}_{\bM_{m_i,n_i}}(A_i)$ whenever $b_{1j}=b_{i1}=0$
for all $i,j \neq 1$.
Similarly, we also have $B=(b_{ij})\in \Span \mathcal{P}_{\bM_{m_i,n_i}}(A_i)$ whenever $b_{2j}=b_{i2}=0$ for all $i,j \neq 2$.
Furthermore, the real symmetric matrix
$\begin{pmatrix} 0&1\\ 1&0\end{pmatrix}\oplus 0_{m_i-2, n_i-2} \in \mathcal{P}_{\bM_{m_i,n_i}}(A_i)$.
It also holds that the complex hermitian matrix
 $\mathrm{i}\begin{pmatrix} 0&1\\ -1&0\end{pmatrix}\oplus 0_{m_i-2, n_i-2} \in \mathcal{P}_{\bM_{m_i,n_i}}(A_i)$ if $\mathbb{F}=\mathbb{C}$. Thus,
$\mathcal{P}_{\bM_{m_i,n_i}}(A_i)$ spans $\bM_{m_i,n_i}$ if $\mathbb{F}=\mathbb{C}$, and $\dim (\Span \mathcal{P}_{\bM_{m_i,n_i}}(A_i))\geq m_i n_i - 1$ if $\mathbb{F}=\mathbb{R}$.
It follows that $\Span \mathcal{P}_{\bM}(A)$ contains the direct summand $\bM_{m_i,n_i}$, except possibly for the one dimensional space spanned by
$\begin{pmatrix} 0&1\\ -1&0\end{pmatrix}\oplus 0_{m_i-2, n_i-2}$ in
the real case.

Fix any nonzero $B_i\in \mathcal{P}_{\bM_{m_i,n_i}}(A_i)$.
For any $j\neq i$, let $B_j\in \bM_{m_j,n_j}$ with $\|B_j\|\leq \|B_i\|$, and let $b_l\in \IF$ with
$|b_l|\leq \|B_i\|$ for $l=1,\ldots, h$.
Then $B=\bigoplus_j B_j\oplus (b_1)\oplus \cdots \oplus (b_l)\in \mathcal{P}_{\bM}(A)$.
In particular, $B'\in \mathcal{P}_{\bM}(A)$, where $B_i'=B_i$, all other $B_j'=0_{m_j, n_j}$
and $b'_l=0$ for $l=1,\ldots,h$.
  Thus $B-B'\in \Span \mathcal{P}_{\bM}(A)$. This shows that all other direct summands
   $\bM_{m_j,n_j}$, as well as the scalar parts $\underbrace{\IF \oplus_\infty \cdots \oplus_\infty \IF }_h$,
   belong to
$\Span\mathcal{P}_{\bM}(A)$ for $j\neq i$. Consequently, the codimension of $\Span\mathcal{P}_{\bM}(A)$ in ${\bM}$ is either $0$ when $\IF=\IC$,
or at most $1$ when $\IF=\IR$.

Suppose now that there are distinct indices $i, j$ such that $s_1(A_i)=s_1(A_j)=\|A\|$. Without loss of generality, we can
assume $i=1, j=2$ and $\|A\|=1$.
Let $U_1, U_2, V_1, V_2$ be unitary matrices of appropriate sizes such that
$U_l^*A_lV_l  = \begin{bmatrix} 1\end{bmatrix}\oplus F_l$
with $\|F_l\|\leq 1$ for $l=1,2$.
Let $\mathrm{E}^{(1)}_{11}\in \bM_{m_1,n_1}$ be the matrix unit with $1$ at its $(1,1)$th entry and zero
elsewhere. Let $B_l\in \bM_{m_l, n_l}$ with $\|B_l\|\leq 1$ for
$l=2,\ldots, k$, and let $b_1,\ldots, b_h\in \IF$ with absolute values not greater than $1$.
Then both $U_1\mathrm{E}^{(1)}_{11}V_1^*\oplus B_2 \oplus B_3 \oplus \cdots \oplus B_k\oplus (b_1)\oplus \cdots\oplus (b_h)$ and
$U_1\mathrm{E}^{(1)}_{11}V_1^*\oplus 0 \oplus 0 \oplus \cdots \oplus 0$ belong to $\mathcal{P}_{\bM}(A)$. Thus their difference $0\oplus B_2 \oplus B_3 \oplus \cdots \oplus B_k\oplus (b_1)\oplus \cdots\oplus (b_h)$ belongs to $\Span \mathcal{P}_{\bM}(A)$.
Consequently, $0\oplus_\infty \bM_{m_2, n_2} \oplus_\infty \bM_{m_3, n_3}  \oplus_\infty \cdots \oplus_\infty \bM_{m_k, n_k}\oplus_\infty
\underbrace{\IF \oplus_\infty \cdots \oplus_\infty \IF }_h\in \Span \mathcal{P}_{\bM}(A)$. In a similar manner, we see that
$\bM_{m_1, n_1}\oplus_\infty 0 \oplus_\infty \bM_{m_3, n_3} \oplus_\infty \cdots \oplus_\infty \bM_{m_k, n_k}\oplus_\infty
\underbrace{\IF \oplus_\infty \cdots \oplus_\infty \IF }_h\in \Span \mathcal{P}_{\bM}(A)$. Hence $\Span \mathcal{P}_{\bM}(A)= {\bM}$ in this case.

If a matrix component and a scalar component both attain $\|A\|$,
the same difference argument, using a peaking rank-one operator in the
matrix component and the scalar coordinate in turn, gives the whole
space. Finally, if two scalar components of $A$ in its diagonal part in $\underbrace{\IF \oplus_\infty \cdots \oplus_\infty \IF }_h$
assume the norm $\|A\|$, then a similar argument as above shows that
$\Span \mathcal{P}_{\bM}(A)= {\bM}$. This finishes the proof of part (1).

Parts (2) and (3) are rather straightforward, and we skip the details.
\end{proof}

We will show that the linear parallelism preserver $T$ preserves
peaking partial isometries of rectangular matrices in some way by using the following
technical lemma.

\begin{lemma}[{\cite[Lemma 2.5]{LTWW-lp}}] \label{lem:D-matrix}
Let $D = (d_{rs}) \in \bM_n$.
Suppose either
 $d_{jj} > |d_{jk}|$ whenever $j \ne k$, or
  $d_{jj} > |d_{kj}|$ whenever $j \ne k$.
The following conditions are equivalent.
\begin{itemize}
\item [{\rm (a)}]
For any
 $\mathbf{x} = (x_1, \dots, x_n)^{\TT}\in \IF^n$ with $\mathbf{y} = D^{\TT}\mathbf{x} = (y_1, \dots, y_n)^{\TT}$, we have
\begin{align*}
\|\mathbf{x}\|_\infty = |x_r| > |x_s|\ \text{whenever  $s \ne r$}\quad
\implies\quad \|\mathbf{y}\|_\infty = |y_r| \geq |y_s|\ \text{whenever $s\ne r$.}
\end{align*}
\item[{\rm (b)}] Either $D = d_{11}I_n$, or $n = 2$ and $D = D^*$ with $d_{11}=d_{22}>|d_{12}|$.
\end{itemize}
\end{lemma}

The final step in proving Theorem \ref{thm:Sum-R} is to verify that
$T$ sends disjoint rank one matrices to disjoint rank one matrices.
Then the following well known result applies.

\begin{lemma}[{\cite[Theorem 7]{BL}}]\label{lem:BL}
  Let $\IF$ be any field containing at least four elements, and $k\in\{1,\ldots,  \min\{m,n\}\}$.
  Let $T: \bM_{m,n}(\IF)\to \bM_{m,n}(\IF)$ be a bijective linear map sending
  rank $k$ rectangular matrices to rank $k$ rectangular matrices. Then there are
  $U\in \bM_m(\IF)$, $V\in \bM_n(\IF)$ such that $T$ assumes either the form
  $A\mapsto UAV$, or the form $A\mapsto UA^\TT V$ which happens only when $m=n$.
\end{lemma}

\begin{proof}[Proof of the implication (a) $\implies$ (d) in Theorem \ref{thm:Sum-R}]
To ease the notations, in the following, we identify the diagonal summand
 $\underbrace{\IF \oplus_\infty \cdots \oplus_\infty \IF }_h$ in $\bM$
 with the finite dimensional supnormed space $\ell_\infty^h$. We
 also write $\mathbf{R}_i= \bM_{m_i, n_i}$ for $i=1,\ldots, k$. Hence,
$$
\bM = \bigoplus_{i=1}^k \mathbf{R}_i \oplus_\infty \ell_\infty^h.
$$

 Replacing $\mathbf{R}_i$ with $\mathbf{R}_i^\TT$ if necessary, we can assume that
$n_i\geq m_i\geq 2$. After a permutation of indices, we can further assume that $m_1+n_1 \geq m_{2}+n_{2} \geq\cdots\geq m_k+n_k\geq 4$,
and if $m_i+n_i = m_{j} + n_{j}$, we have
$m_i\geq m_{j}$ whenever $i<j$. For any element $A\in\bM$, we write
$A= A_{{1}}\oplus \cdots \oplus A_{{k}}\oplus  (a_1,\ldots, a_h)$ with
the block components $A_{{i}}\in \mathbf{R}_i$
for $i=1, \ldots, k$, and $(a_1,\ldots, a_h)\in \ell_\infty^h$.

By assumption,
$$
T(\mathcal{P}_{\bM}(A)) \subseteq \mathcal{P}_{\bM}(TA),
$$
and thus
$$
T(\Span {\mathcal{P}_{\bM}}(A))  \subseteq \Span {\mathcal{P}_{\bM}}(TA) \quad\text{for every $A\in {\bM}$}.
$$
Since $T$ is bijective, the codimensions of the subspaces spanned by
$\mathcal{P}_{\bM}(A)$ and $\mathcal{P}_{\bM}(TA)$ of $\bM$ satisfy
\begin{align}\label{eq:codim}
 \operatorname{codim} \Span\mathcal{P}_{\bM}(TA) \leq \operatorname{codim}
\Span \mathcal{P}_{\bM}(A)\quad\text{for all $A\in {\bM}$.}
\end{align}

\noindent\textbf{Step 1.}
We start with the block component $\mathbf{R}_1=\bM_{m_1,n_1}$.
Let $\mathrm{E}^{(1)}_{ij}\in \mathbf{R}_1$ be the matrix unit with $1$ at its $(i,j)$th entry
 and zero elsewhere for $1\leq i\leq m_1$ and $1\leq j\leq n_1$; we also regard
 $\mathrm{E}^{(1)}_{ij}=\mathrm{E}_{ij}\in \mathbf{M}_{m+h, n+h}$ as an element in $\mathbf{M}$.
 If
 $F_{11}=T^{-1}(\mathrm{E}^{(1)}_{11})$ then
 $T(\Span {\mathcal{P}_{\bM}}(F_{11}))  \subseteq \Span {\mathcal{P}_{\bM}}(\mathrm{E}^{(1)}_{11})$.
 In view of Lemma \ref{lem:Pk1}(2a), $\Span \mathcal{P}_\bM(\mathrm{E}^{(1)}_{11})$
has codimension $m_1+n_1-2\geq 2$ in $\bM$. It follows from Lemma \ref{lem:Pk1}(3) that $F_{11}$ has peak rank one and
 a peaking partial isometry of $F_{11}$
 belongs to a block
$\mathbf{R}_s=\bM_{m_s,n_s}$.
By \eqref{eq:codim},
   the codimension of $\Span {\mathcal{P}_{\bM}}(F_{11})$
in ${\bM}$ equals $m_s + n_s - 2
\geq m_1 + n_1 - 2$.
  Because $m_1 + n_1$ is the biggest among others by construction,
we see that
$$
m_s + n_s = m_1 + n_1\quad\text{as well as}\quad m_s\leq m_1.
$$
   Consequently,
\begin{align}\label{eq:TSPF11}
T(\Span {\mathcal{P}_{\bM}}(F_{11})) =\Span {\mathcal{P}_{\bM}}(\mathrm{E}^{(1)}_{11}).
\end{align}

Consider $\mathrm{E}^{(1)}_{12}$ in $\mathbf{R}_1$.
Argue as in above, $F_{12}:=T^{-1}(\mathrm{E}^{(1)}_{12})$ has peak rank one. If a
peaking partial isometry of $F_{12}$
belongs to any block component other than that for $F_{11}$,
we have $F_{11} + \alpha F_{12}=T^{-1}(\mathrm{E}^{(1)}_{11}+\alpha \mathrm{E}^{(1)}_{12})$
has peak rank two for some scalar $\alpha$.
This is a contradiction,
since the rank one matrix
$\mathrm{E}^{(1)}_{11}+\alpha \mathrm{E}^{(1)}_{12}$ always has peak rank one, and
thus $T^{-1}(\mathrm{E}^{(1)}_{11}+\alpha \mathrm{E}^{(1)}_{12})$ cannot have peak rank greater than one
by \eqref{eq:codim} and Lemma \ref{lem:Pk1}(1).
It thus ensures that
the peaking partial isometries of $T^{-1}(\mathrm{E}^{(1)}_{11})$ and $T^{-1}(\mathrm{E}^{(1)}_{12})$
belong to the same block component $\mathbf{R}_s$.
Applying the same reasoning to the pair $\mathrm{E}^{(1)}_{22}$ and $\mathrm{E}^{(1)}_{12}$, we
see that all the peaking partial isometries of $T^{-1}(\mathrm{E}^{(1)}_{11})$,
$T^{-1}(\mathrm{E}^{(1)}_{12})$ and $T^{-1}(\mathrm{E}^{(1)}_{22})$
belong to the same block component $\mathbf{R}_s$. Eventually, for $\mathrm{E}^{(1)}_{ij}$ in $\mathbf{R}_1$,
its preimage $F_{ij}=T^{-1}(\mathrm{E}^{(1)}_{ij})\in \operatorname{Pr}_1$ has a peaking partial isometry
in the block component $\mathbf{R}_s$.

By Lemma \ref{lem:Pk1}(2a),
there are unitary matrices $U_s\in\bM_{m_s}$ and $V_s\in\bM_{n_s}$ such that
$$
F_{11}=\bigoplus_i F_{11,i} \oplus (f_{11,1},\ldots, f_{11,h})\quad\text{with}\quad
 F_{11,s} = U_s([d_{11}] \oplus G_1)V_s,
 $$
in which
 $$
 d_{11}> \|G_1\|, \|F_{11,i}\|, |f_{11,j}| \quad\text{for all $i\neq s$ and $j=1,\ldots, h$.}
   $$
   It follows
\begin{align}\label{eq:SP11}
\hspace{1cm} &\Span \mathcal{P}_\bM(F_{11})=  \Bigg\{\bigoplus_i B_{i}\oplus (b_1) \oplus\cdots\oplus (b_h)\in \bM: B_{s} = U_s ([b] \oplus F_1)V_s \notag\\
&\hspace{1cm} \text{with arbitrary}\ B_{i}\in\mathbf{R}_i \ \text{for any}\ i\neq s,
 F_1 \in \bM_{m_s-1,n_s-1}\ \text{and}\ b, b_1,\ldots, b_h\in \IF\Bigg\}.
\end{align}

Since $\mathrm{E}^{(1)}_{22} \in \Span \mathcal{P}_\bM(\mathrm{E}^{(1)}_{11})$,
  the invertibility of $T$ and \eqref{eq:TSPF11} imply that
$F_{22}=T^{-1}(\mathrm{E}^{(1)}_{22}) \in \Span \mathcal{P}_\bM(F_{11})$. By \eqref{eq:SP11}, the $s$th component of $F_{22}$
 has the form
 $$
 F_{22,s}=U_s([d_{21}] \oplus C_1)V_s\in \mathbf{R}_s\quad\text{with}\quad C_1\in
 \bM_{m_s-1,n_s-1}.
  $$
  Since $T(F_{11}) = \mathrm{E}^{(1)}_{11}$ and $T(F_{22}) = \mathrm{E}^{(1)}_{22}$
are not parallel, $F_{11}$ and $F_{22}$ are not parallel.
Because $F_{22}\in \operatorname{Pr}_1$ has a peaking partial isometry in $\mathbf{R}_s$,
we have
$$
\|C_1\|=\|F_{22}\| > |d_{21}|\ \text{and}\ \|C_1\|> \|F_{22,i}\|, |f_{22,j}|\quad
\text{for $i\neq s$ and $j=1,\ldots, h$.}
$$
 We may change $U_s, V_s$ to
$U_s([1]\oplus U')$ and $([1]\oplus V')V_s$ for suitable unitary matrices $U' \in \bM_{m_s-1}$ and
$V' \in \bM_{n_s-1}$
to arrive at that
$$
F_{11,s} = U_s([d_{11}] \oplus G_1')V_s\quad
\text{and}\quad
F_{22,s} = U_s([d_{21}] \oplus [d_{22}] \oplus C_2)V_s
$$
where
$d_{11}=\|F_{11}\|>\|G_1'\|$, $d_{22} =\|F_{22}\| > |d_{21}|$ and $d_{22} > \|C_2\|$.

Since $\mathrm{E}^{(1)}_{11}\in \Span\mathcal{P}_\bM(\mathrm{E}^{(1)}_{22})$, we see that
$F_{11}\in \Span \mathcal{P}_\bM(F_{22})$ but not parallel with $F_{22}$.
Observe that $B\in \Span \mathcal{P}_\bM(F_{22})$ exactly when the $s$th component
$B_{s}$ of $B$ has zero at the entries in the second row and the second column
other than the $(2,2)$th one.
We can then refine the above observation to
\begin{align*}
F_{11, s} &= U_s([d_{11}] \oplus [d_{12}] \oplus   G_2)V_s\quad\text{with}
\quad \| G_2\|< \|F_{11}\|=d_{11}\quad\text{and}\quad |d_{12}|< d_{11}, \\
F_{22, s} &= U_s([d_{21}] \oplus [d_{22}] \oplus   C_2)V_s\quad\text{with}
\quad \| C_2\|< \|F_{22}\|=d_{22}\quad\text{and}\quad |d_{21}|< d_{22}.
\end{align*}

Continue this line of reasoning with $\mathrm{E}^{(1)}_{33}, \ldots, \mathrm{E}^{(1)}_{m_1 m_1}$.
We can assert that
$F_{ii} = T^{-1}(\mathrm{E}^{(1)}_{ii})\in \operatorname{Pr}_1$ with peaking partial isometries in $\mathbf{R}_s$,
and there are unitary matrices $U_s\in \bM_{m_s}$ and $V_s\in\bM_{n_s}$ such that
the components of $F_{ii}$ in $\mathbf{R}_s$ assume the form
\begin{align*}
F_{ii, s}  &= U_s([d_{i1}] \oplus \cdots\oplus [d_{im_1}]\oplus \tilde{F}_i)V_s \notag\\
&\hspace{1cm}\text{with}\quad
\|\tilde{F}_i\|<\|F_{ii}\|=d_{ii}\quad\text{and}\quad |d_{ij}|<d_{ii} \quad\text{for all $j\neq i$ and $i=1,2,\ldots, m_1$.}
\end{align*}
In particular, we see that $m_1\leq m_s$.
By the assumption on the indices, we have $m_1=m_s$, and thus $(m_s, n_s)=(m_1,n_1)$
 since $m_s+n_s=m_1+n_1$.
After a permutation of indices, we can assume that $s=1$.

At this point, we can write
\begin{align*}
T^{-1}(\mathrm{E}^{(1)}_{ii})  &=  U_1 \begin{pmatrix} \begin{matrix} d_{i1} &&&  \\
&d_{i2} &&  \\ &&\ddots & \\ &&&   d_{im_1}\end{matrix} &
\multicolumn{3}{c}{\raisebox{0ex}[0pt]{\Huge0}_{m_1, n_1-m_1}}\end{pmatrix}V_1
\oplus F_{ii,2}\oplus\cdots\oplus F_{ii,k}\oplus (f_{ii,1}, \ldots, f_{ii,h})
\notag\\[10pt]
&\hspace{2cm}    \text{in which}\quad
 \|F_{ii,r}\|, |f_{ii,s}|,   |d_{ij}|<d_{ii}
 \quad  \text{for all $r\geq 2$, $s\geq 1$,
and $j\neq i$.}
\end{align*}
Moreover,
\begin{align}\label{eq:TSP(1)}
 \Span {\mathcal{P}_{\bM}}(T^{-1}(\mathrm{E}^{(1)}_{ii}))  =T^{-1}\left(\Span {\mathcal{P}_{\bM}}(\mathrm{E}^{(1)}_{ii})
 \right) \quad\text{for $i=1,\ldots, m_1$,}
\end{align}
and
\begin{align}\label{eq:SP(1)}
 \hspace{1cm} \bigcap_{i=1}^{m_1}\, & \Span \mathcal{P}_\bM(T^{-1}(\mathrm{E}^{(1)}_{ii})) =  \notag \\
&\hspace{0.5cm}\Biggl\{ U_1\begin{pmatrix} \begin{matrix} b_{1 1} &&&  \\
& b_{12} &&  \\ &&\ddots & \\ &&&    b_{1   m_1}\end{matrix} &
\multicolumn{3}{c}{\raisebox{0ex}[0pt]{\Huge0}_{m_1, n_1-m_1}}\end{pmatrix}V_1
\oplus B_2\oplus\cdots\oplus B_k\oplus  (b_{1}, \ldots, b_{h}) \notag \\
&\hspace{0.75cm}  \text{with arbitrary}\
 B_j  \in \bM_{m_j,n_j}\ \text{for $j\neq 1$, and}\
 b_{1 1}, b_{1 2}, \ldots, b_{1 m_1}, b_1,\ldots, b_h\in \IF\Biggr\}.
\end{align}

\noindent\textbf{Step 2.}
Consider the matrix unit $\mathrm{E}^{(2)}_{ij}$ in $\mathbf{R}_2=\bM_{m_2,n_2}$
 with the $(i,j)$th entry equal $1$ and others equal $0$ for $1\leq i\leq m_2$ and $1\leq j\leq n_2$; we also regard
 $\mathrm{E}^{(2)}_{ij}=\mathrm{E}_{m_1+i, n_1+j}\in \mathbf{M}_{m+h, n+h}$ as an element in $\mathbf{M}$.
Since $\mathrm{E}^{(2)}_{ij}\in \Span \mathcal{P}_\bM(\mathrm{E}^{(1)}_{ll})$,
 we see that
$T^{-1}(\mathrm{E}^{(2)}_{ij})\in \Span \mathcal{P}_\bM(T^{-1}(\mathrm{E}^{(1)}_{ll}))$ for all feasible
indices $i,j, l$
by \eqref{eq:TSP(1)}.
It follows from \eqref{eq:SP(1)} that
 $$
 T^{-1}(\mathrm{E}^{(2)}_{ij})= U_1(\sum_{l=1}^{m_1} d^{(2,1)}_{ij,l} \mathrm{E}^{(1)}_{ll})V_1\oplus F^{(2,2)}_{ij} \oplus \cdots\oplus F^{(2,k)}_{ij}\oplus (f^{(2)}_{ij,1}, \ldots, f^{(2)}_{ij,h})
 $$
 for some scalars $d^{(2,1)}_{ij,1}$, \ldots,  $d^{(2,1)}_{ij,m_1}$, $f^{(2)}_{ij,1}$,
\ldots, $f^{(2)}_{ij,h}$, and
 $F^{(2,p)}_{ij}\in \bM_{m_p,n_p}$, $p=2,\ldots, k$.
Since $\mathrm{E}^{(2)}_{ij}$ is not parallel with
 $\mathrm{E}^{(1)}_{ll}$, we see that $T^{-1}(\mathrm{E}^{(2)}_{ij})$ is not parallel with $T^{-1}(\mathrm{E}^{(1)}_{ll})$.
Together with Lemma \ref{lem:Pk1}, we see that $T^{-1}(\mathrm{E}^{(2)}_{ij})$ has peak rank one, and
there is a $q_{ij}\geq 2$ such that $T^{-1}(\mathrm{E}^{(2)}_{ij})$ has a peaking partial isometry
in $\mathbf{R}_{q_{ij}}$ such that $|d^{(2,1)}_{ij,l}|< \|F^{(2,q_{ij})}_{ij}\|=
 \|T^{-1}(\mathrm{E}^{(2)}_{ij})\|$ for $l=1, \ldots, m_1$, $\|F^{(2,t)}_{ij}\|< \|F^{(2,q_{ij})}_{ij}\|$
  for $t\neq q_{ij}$, and $|f^{(2)}_{ij,r}|< \|F^{(2,q_{ij})}_{ij}\|$ for $r=
 1,\ldots, h$. Due to \eqref{eq:codim}, such $q_{ij}$
 satisfy $m_{q_{ij}} + n_{q_{ij}} -2 \geq m_2 + n_2-2$. Because $m_2+n_2$ is the biggest one among all $m_r + n_r$ for $r=2,\ldots, k$, we see that $m_{q_{ij}} + n_{q_{ij}} = m_2 + n_2$,
 as well as $m_{q_{ij}}\leq m_2$.
 Arguing as in Step 1, we see that all indices
 $q_{ij}$ are the same, $q$ say, such that $m_q=m_2$ and $n_q=n_2$.
 We can assume that $q=2$,
after a permutation of the direct summands of $\bM$ fixing $\mathbf{R}_1$,
and there are
  unitary matrices $U_2\in \bM_{m_2}$ and $V_2\in \bM_{n_2}$ such that
\begin{align*}
 T^{-1}(\mathrm E^{(2)}_{ii})
 &=\bigoplus_{r=1}^{2}U_r\left(\sum_{t=1}^{m_r}
 d^{(2,r)}_{it}\mathrm E^{(r)}_{tt}\right)V_r
 \oplus F^{(2,3)}_{ii}\oplus\cdots\oplus F^{(2,k)}_{ii}
 \oplus(f^{(2)}_{ii,1},\ldots,f^{(2)}_{ii,h}).
\end{align*}
\text{for $i=1,\ldots, m_2$,} in which
 $$
 \|T^{-1}(\mathrm{E}^{(2)}_{ii})\|=
 d^{(2,2)}_{ii}> |d^{(2,2)}_{ij}|,\ |d^{(2,1)}_{il}|,\ \|F_{ii}^{(2,t)}\|,\ |f^{(2)}_{ii,r}|
$$
 \text{for
 $j\neq i$, $l=1,\ldots, m_1$, $t= 3, \ldots, k$, and $r=1,\ldots, h$.}

Moreover,
\begin{align*}
 \Span {\mathcal{P}_{\bM}}(T^{-1}(\mathrm{E}^{(2)}_{ii}))  =T^{-1}\left(\Span {\mathcal{P}_{\bM}}(\mathrm{E}^{(2)}_{ii})
 \right) \quad\text{for $i=1,\ldots, m_2$,}
\end{align*}
and
\begin{align*}
 &\bigcap_{r=1}^{2}\bigcap_{i=1}^{m_r}
 \Span\mathcal P_\bM(T^{-1}(\mathrm E^{(r)}_{ii}))\\
 &\hspace{1cm}=\Biggl\{\bigoplus_{r=1}^{2}U_r\left(\sum_{t=1}^{m_r}
 b^{(2,r)}_t\mathrm E^{(r)}_{tt}\right)V_r
 \oplus B_3\oplus\cdots\oplus B_k\oplus(b_1,\ldots,b_h): \\
 &\hspace{3cm}  B_j\in\bM_{m_j,n_j}\ (j\geq3),\quad
 b^{(2,r)}_t,b_l\in\IF\Biggr\}.
\end{align*}

Since $\mathrm{E}^{(1)}_{ii}\in \Span\mathcal{P}_\bM(\mathrm{E}^{(2)}_{ll})$,
we see that $T^{-1}(\mathrm{E}^{(1)}_{ii})\in \Span\mathcal{P}_\bM(T^{-1}(\mathrm{E}^{(2)}_{ll})$
for all feasible indices $i, l$. On the other hand,
 $\mathrm{E}^{(1)}_{ii}$ and $\mathrm{E}^{(2)}_{ll}$ are not parallel,
so are not $T^{-1}(\mathrm{E}^{(1)}_{ii})$ and $T^{-1}(\mathrm{E}^{(2)}_{ll})$ for $i=1,\ldots, m_1$ and
$l=1,\ldots, m_2$. As in the previous step, we can refine the expressions as
\begin{align*}
 T^{-1}(\mathrm E^{(r)}_{ii})
 &=\bigoplus_{s=1}^{2}U_s\left(\sum_{t=1}^{m_s}
 d^{(r,s)}_{it}\mathrm E^{(s)}_{tt}\right)V_s\oplus F^{(r,3)}_{ii}\oplus\cdots\oplus F^{(r,k)}_{ii}
 \oplus(f^{(r)}_{ii,1},\ldots,f^{(r)}_{ii,h})
\end{align*}
 for $r=1,2$ and $i=1,\ldots,m_r$; in which
\begin{align*}
 d^{(1,1)}_{ii} &> |d^{(1,s)}_{it}|, \|F^{(1,r)}_{ii}\|, |f^{(1)}_{ii,u}|\quad\text{whenever
 $(s, t)\neq (1, i)$, $r\geq 3$, $u=1,\ldots,h$, and} \\
  d^{(2,2)}_{jj} &> |d^{(2,s)}_{jt}|, \|F^{(2,r)}_{jj}\|, |f^{(2)}_{jj,u}|\quad\text{
  whenever $(s, t)\neq (2, j)$, $r\geq 3$, $u=1,\ldots,h$}.
 \end{align*}

\noindent\textbf{Step 3.}
Continue this process, after possibly permutations on indices,
we shall eventually arrive at
\begin{multline}\label{eq:Tii-form}
 T^{-1}(\mathrm{E}^{(r)}_{ii})   = \\
  U_1\left(\sum_{t=1}^{m_1} d^{(r,1)}_{it} \mathrm{E}^{(1)}_{tt}\right)V_1\oplus
  U_2\left(\sum_{t=1}^{m_2} d^{(r,2)}_{it} \mathrm{E}^{(2)}_{tt}\right)V_2 \oplus \cdots\oplus
  U_k\left(\sum_{t=1}^{m_k} d^{(r,k)}_{it} \mathrm{E}^{(k)}_{tt}\right)V_k\oplus (f^{(r)}_{ii,1}, \ldots, f^{(r)}_{ii,h}),
 \end{multline}
 for $r=1,\ldots, k$; in which $U_i\in\bM_{m_i}$ and $V_i\in\bM_{n_i}$ are unitary matrices, and
\begin{align*}
 d^{(r,r)}_{ii}  > |d^{(r,s)}_{it}|,   |f^{(r)}_{ii,u}|\quad\text{whenever
 $(s,t)\neq (r,i)$  and $u=1,\ldots,h$.}
 \end{align*}

\noindent\textbf{Step 4.}
We study the diagonal summands $\mathbf{e}_{\alpha}=
(0, \ldots, 0, \underbrace{1}_{\alpha\text{th place}}, 0, \ldots, 0)$ in
$\ell_\infty^h$ for $\alpha=1,\ldots, h$.
Since $\mathbf{e}_{\alpha}\in \Span \mathcal{P}_\bM(\mathrm{E}^{(r)}_{ii})$ but not
parallel with $\mathrm{E}^{(r)}_{ii}$ for any $\alpha=1,\ldots, h$, $r=1, \ldots, k$ and
 $i=1,\ldots, m_r$,
we see that $T^{-1}(\mathbf{e}_{\alpha})$ assumes the form \eqref{eq:Tii-form} but
attains its norm at some components in $\ell_\infty^h$.
We claim that $T^{-1}(\mathbf{e}_{\alpha})$ has peak rank one.
Indeed, since each of the $h$ many
block matrices $T^{-1}(\mathbf{e}_1), \ldots,  T^{-1}(\mathbf{e}_h)$
attain their norms at some of the $h$ many
coordinates from $\ell_\infty^h$, if any of them attains its norm at more than one coordinates,
there will be distinct $T^{-1}(\mathbf{e}_{\alpha})$ and
$T^{-1}(\mathbf{e}_{\beta})$
  attain their norms at a same coordinate, and thus they are parallel
  with each other. Since $T$ preserves parallel pairs,
we have a contradiction that $\mathbf{e}_{\alpha}$ and $\mathbf{e}_{\beta}$ are parallel.
After a possibly permutation of the coordinates in the direct summand
  $\ell_\infty^h$, we can find unimodular scalars $e^{\mathrm{i}\theta_{1}}, \ldots, e^{\mathrm{i}\theta_{h}}$ such that
\begin{multline}\label{eq:Teu}
 T^{-1}(\mathbf{e}_{\alpha})  = U_1\left(\sum_{t=1}^{m_1} d^{(k+\alpha,1)}_{tt} \mathrm{E}^{(1)}_{tt}\right)V_1
 \oplus
  U_2\left(\sum_{t=1}^{m_2} d^{(k+\alpha,2)}_{tt} \mathrm{E}^{(2)}_{tt}\right)V_2
   \oplus \cdots\oplus
  U_k\left(\sum_{t=1}^{m_k} d^{(k+\alpha,k)}_{tt} \mathrm{E}^{(k)}_{tt}\right)V_k \\
   \oplus (e^{\mathrm{i}\theta_{1}}f^{(k+\alpha)}_1, \ldots, e^{\mathrm{i}\theta_{h}}f^{(k+\alpha)}_h),
\end{multline}
 for $\alpha=1,\ldots, h$; in which
\begin{align*}
f^{(k+\alpha)}_{\alpha}  > |d^{(k+\alpha, s)}_{tt}|,\ |f^{(k+\alpha)}_{l}| \quad\text{for all feasible indices $s,t$, and $l\neq \alpha$.}
 \end{align*}

\noindent\textbf{Step 5.}
We now rename the disjoint
norm one partial isometries
$\mathrm{E}^{(1)}_{11}, \ldots,  \mathrm{E}^{(1)}_{m_1m_1}$ in $\mathbf{R}_1$,
$\mathrm{E}^{(2)}_{11}$, \ldots,  $\mathrm{E}^{(2)}_{m_2m_2}$ in $\mathbf{R}_2$, \ldots,
$\mathrm{E}^{(k)}_{11}, \ldots,  \mathrm{E}^{(k)}_{m_km_k}$ in $\mathbf{R}_k$, and
$\mathbf{e}_1, \ldots, \mathbf{e}_h$
in $\ell_\infty^h$, as
$$
w_1, \ldots, w_{m_1}, w_{m_1+1}, \ldots, w_{m_1+m_2}, \ldots, w_{m}, w_{m+1}, \ldots, w_{m+h},
$$
where $m=m_1+m_2 + \cdots + m_k$. Correspondingly, we rename the disjoint
norm one partial isometries
$U_1\mathrm{E}^{(1)}_{11}V_1, \ldots, U_1\mathrm{E}^{(1)}_{m_1m_1}V_1$ in $\mathbf{R}_1$,
$U_2\mathrm{E}^{(2)}_{11}V_2, \ldots, U_2\mathrm{E}^{(2)}_{m_2m_2}V_2$ in $\mathbf{R}_2$, \ldots,
$U_k\mathrm{E}^{(k)}_{11}V_k, \ldots, U_k\mathrm{E}^{(k)}_{m_km_k}V_k$ in $\mathbf{R}_k$, and $e^{\mathrm{i}\theta_{1}}\mathbf{e}_1, \ldots, e^{\mathrm{i}\theta_{h}}\mathbf{e}_h$
in $\ell_\infty^h$, as
$$
w'_1, \ldots, w'_{m_1}, w'_{m_1+1}, \ldots, w'_{m_1+m_2}, \ldots, w'_{m}, w'_{m+1}, \ldots, w'_{m+h}.
$$
Summarizing \eqref{eq:Tii-form} and \eqref{eq:Teu},
as well as renaming $d^{(r,s)}_{ij}$'s and $f^{(k+\alpha)}_t$'s as $d_{ij}$'s
(with $d_{ij}=0$ for infeasible/missing indices), we can write
$$
T^{-1}(w_i)= \sum_{j=1}^{m+h} d_{ij}w'_j\quad\text{for $i=1,\ldots, m+h$},
$$
in which
$$
d_{ii} > |d_{ij}| \quad\text{whenever $i\neq j$.}
$$

For any scalar vector $\mathbf{x}=(x_1, \ldots, x_{m+h})^\TT$, we associate it with
the diagonal matrix
$$
A(\mathbf{x}) = \sum_{i=1}^{m+h} x_i w_i\in \bM.
$$
It is clear that $\|A(\mathbf{x})\|=\|x\| = \max_i |x_i|$.
If $\|x\|=|x_i|>|x_j|$ whenever $i\neq j$
 then $A(\mathbf{x})$ is parallel with the matrix $w_i$ but not any other $w_j$'s.
Observe that
\begin{align*}
T^{-1}(A(\mathbf{x})) &= \sum_{i=1}^{m+h} x_i T^{-1}(w_i) = \sum_{i=1}^{m+h} \left(x_i \sum_{j=1}^{m+h} d_{ij}w'_j\right)\\
&= \sum_{j=1}^{m+h}\left(\sum_{i=1}^{m+h} x_id_{ij}\right)w'_j = \sum_{j=1}^{m+h} y_j w'_j,
\end{align*}
where $y_j = \sum_{i=1}^{m+h} x_id_{ij}$, or
$$
\mathbf{y} = (y_1, \ldots, y_{m+h})^\TT = D^\TT \mathbf{x}\quad \text{with $D=(d_{ij})\in\bM_{m+h}$.}
$$

If $j\ne i$, preservation of parallelism by $T$ implies that
$T^{-1}(A(\mathbf x))$ is not parallel with $T^{-1}(w_j)$.
The latter has a unique maximal coordinate at $w'_j$, by
\eqref{eq:Tii-form} and \eqref{eq:Teu}. Hence the diagonal operator
$T^{-1}(A(\mathbf x))=\sum_jy_jw'_j$ cannot attain its norm at $w'_j$
for any $j\ne i$.
Consequently,
 $T^{-1}(A(\mathbf{x}))$ is parallel with only
 $w'_i$, and thus $\|\mathbf{y}\|=|y_i|> |y_j|$ whenever $i\neq j$. In summary, $\mathbf{y} = D^\TT \mathbf{x}$ satisfies the condition
$$
\|\mathbf{x}\|=|x_i|> |x_j|\ \text{whenever $i\neq j$}\quad\implies\quad
\|\mathbf{y}\|=|y_i|> |y_j|\ \text{whenever $i\neq j$}.
$$
Since $m+h\geq 3$, by Lemma \ref{lem:D-matrix} we see that $D=d_{11}I_{m+h}$ is a
positive scalar multiple of the
$(m+h)\times (m+h)$ identity
matrix. In other words,
\begin{align}
T^{-1}(\mathrm{E}^{(r)}_{ii}) &= d_{11}U_r \mathrm{E}^{(r)}_{ii}V_r\in \mathbf{R}_r\quad\text{for all $r=1,\ldots,k$, $i=1,\ldots,m_r$,} \label{eq:T-1ii}
\intertext{and}
T^{-1}(\mathbf{e}_{\alpha}) &= d_{11}e^{\mathrm{i}\theta_\alpha}\mathbf{e}_{\alpha} \quad\text{for
$\mathbf{e}_{\alpha}\in \ell_\infty^h$ with $\alpha=1,\ldots, h$.} \label{eq:T-1eu}
\end{align}

\noindent\textbf{Step 6.}
With \eqref{eq:T-1eu} we see that $T^{-1}$ sends the direct summand $\ell_\infty^h$ into $\ell_\infty^h$.
The induced map $\Psi:\ell_\infty^h\to \ell_\infty^h$ obtained
from restricting $T$ to $\ell_\infty^h$ clearly maps $\ell_\infty^h$ onto $\ell_\infty^h$
and, after restoring the permutation on indices, assumes the form
$$
 (a_1, \ldots, a_h)\mapsto \gamma  (u_1a_{\tau(1)}, \ldots, u_h a_{\tau(h)}),
$$
where $\tau$ is a permutation on $\{1, \ldots, h\}$, $\gamma =d_{11}^{-1}>0$,
and $u_\alpha = e^{-\mathrm{i}\theta_\alpha}$ for $\alpha=1, \ldots, h$.

Repeat Steps 1--5 after arbitrary unitary changes of bases in the
matrix blocks. Every rank-one partial isometry can be included in such a
blockwise frame, and every disjoint pair can be included in a common frame.
It follows that $T^{-1}$ sends rank-one matrices in a matrix block to
rank-one matrices in some matrix block, and sends disjoint pairs to
disjoint pairs.

For a fixed $r$, the set of nonzero rank-one matrices in $\mathbf R_r$
is connected, over either $\IR$ or $\IC$, since $m_r,n_r\geq2$.
Its image under the continuous finite-dimensional map $T^{-1}$ cannot
switch between distinct blocks without passing through zero.
Thus the target block is fixed. Equation~\eqref{eq:T-1ii} identifies it
as $\mathbf R_r$. By linearity and dimension,
$T^{-1}(\mathbf R_r)=\mathbf R_r$.

Apply Lemma~\ref{lem:BL} to this restriction of $T^{-1}$. It has the
form $A\mapsto PAQ$ or, for a square block, $A\mapsto PA^{\TT}Q$,
with $P,Q$ invertible. Preservation of disjoint rank-one pairs implies
that $P$ and $Q^*$ preserve orthogonality of vectors: fix orthogonal
vectors in one factor and vary the orthogonal pair in the other.
An orthogonality-preserving linear map on a Hilbert space is a scalar
multiple of an isometry, so $P,Q$ are scalar multiples of unitaries.
Equation~\eqref{eq:T-1ii} fixes their product scale as
$d_{11}=\gamma^{-1}$, independently of the block.
Consequently the restriction of $T$ to $\mathbf R_r$ has the form
$A\mapsto\gamma U_rAV_r$ or $A\mapsto\gamma U_rA^{\TT}V_r$.

Restoring the permutations and unitary similarities, we see that $T$ is the direct sum of the maps
$\Phi_1: \mathbf{R}_{\sigma(1)}\to \mathbf{R}_1$, \ldots,
$\Phi_k: \mathbf{R}_{\sigma(k)}\to \mathbf{R}_k$, and $\Psi:\ell_\infty^h\to \ell_\infty^h$,
which assume the stated forms.
The proof is complete.
\end{proof}

\section[Standard operator spaces]{Linear parallelism preservers of standard operator spaces and \texorpdfstring{$B(H,K)$}{B(H,K)}}\label{s:B(H,K)}

In this section, we consider linear parallelism/TEA preservers of $B(H,K)$,
the Banach space of bounded linear operators between Hilbert spaces,
in which the (real or complex) Hilbert spaces $H, K$ can have arbitrary dimensions.

Let $A^\TT, A^*\in B(K,H)$ denote the transpose and the adjoint of an operator $A\in B(H,K)$
with respect to some arbitrary but fixed orthonormal bases $\{e_i\}_i, \{f_j\}_j$
 of $H$ and $K$; namely, $\langle A^\TT f_j,e_i \rangle_H = \langle A e_i, f_j\rangle_K$
 and $\langle A^* f_j, e_i\rangle_H = \langle f_j , A e_i\rangle_K$
 for all indices $i,j$.
When the underlying field $\IF=\IR$, we have $A^\TT = A^*$; in this case,
  by unitary matrices/operators we refer to real orthogonal matrices/operators. Two operators $A,B\in B(H,K)$ are said to be disjoint if
  $AB^*=0$ and $A^*B=0$; in other words, $A, B$ have orthogonal initial spaces and orthogonal range spaces.

For $u\in K$ and $v\in H$, denote by $u\otimes v$ the rank at most one operator defined by
$(u\otimes v)(\cdot) = \langle \cdot,v\rangle u$.
When $H,K$ have finite dimensions $n$ and $m$, respectively, the matrix representation of $u\otimes v$
is the $m\times n$ rectangular matrix $uv^*$. Note that $u\otimes v$ is a partial isometry, i.e.,
$(u\otimes v)(u\otimes v)^*(u\otimes v)= u\otimes v$, exactly
when $u\otimes v=0$ or $\|u\|\|v\|=1$. In particular, all matrix units $E_{ij}= f_i\otimes e_j$ are partial isometries, where
$\{e_i\}_i$ and $\{f_j\}_j$ are (not necessarily countable) orthonormal bases of $H$ and $K$, respectively.
To ease the notations, we write
$E_{11}, E_{23}, \ldots$ to indicate some `matrix units' $f_1\otimes e_1$, $f_2\otimes e_3$, \ldots,
without really referring to the orders of the (possibly uncountable) orthonormal bases.

Let $A \in B(H,K)$. Recall the \emph{$s$-numbers}
$$
s_1(A) = \|A\|\quad\text{and}\quad  s_k(A) = \inf\{\|A-B\|: \operatorname{rank}(B) < k\}, \quad
\forall k = 2, 3, \ldots.
$$
If $A$ has finite rank then $s_k(A)$ agrees with the $k$th singular value of $A$.
In general, the condition $s_1(A) > s_2(A)$ is equivalent to the existence of a pair of unit vectors
$v\in H$ and $u\in K$ such that
$$
A = s_1(A) u\otimes v + A',
$$
 where $A'$ is
disjoint from the rank one partial
isometry $u\otimes v$, defined by $(u\otimes v)x= \langle x, v\rangle u$, with
$\|A'\|=s_2(A) < s_1(A) = \|A\|$.
In this case, we have $\|Av\|=\|A^*u\| =\|A\|$.
Fix arbitrary unit vectors $e_1\in H$ and $f_1\in K$, and extend them to
(not necessarily countable) orthonormal bases of $H$ and $K$, respectively.
Let $U\in B(K)$ and $V\in B(H)$ be unitary operators such that $U^*f_1=u$ and $Ve_1=v$.
By abusing notations, we can write
$$
UAV = U(s_1(A) u\otimes v + A')V=
 s_1(A)\big( f_1\otimes e_1 + A_2\big) = s_1(A)\big( [1] \oplus A_2\big)
$$
for some operator $A_2$ disjoint from $E_{11}=f_1\otimes e_1$
with $\|A_2\|< 1$.

In general, $A\in B(H,K)$ has {finite} peak rank $p<\infty$ if
$s_1(A)=s_2(A)=\cdots = s_p(A)> s_{p+1}(A)$. In this case,
there are orthogonal unit vectors $v_1, \ldots, v_p$ in $H$ and $u_1, \ldots, u_p$ in $K$ such that
$$
A = s_1(A) \bigl(u_1\otimes v_1 + \cdots + u_p\otimes v_p\bigr) + A'',
$$
where the operator $A''$ is disjoint from all $u_1\otimes v_1$, \ldots, $u_p\otimes v_p$ with
$\|A''\|= s_{p+1}(A) < s_p(A)=\|A\|$.
In a similar fashion, for any orthonormal sets $\{e_1, e_2, \ldots, e_p\}$ in $H$ and $\{f_1, f_2, \ldots, f_p\}$ in $K$, there are unitaries $U\in B(K)$
and $V\in B(H)$ such that
\begin{align}\label{eq:UAV}
UAV = U\left(s_1(A) \bigl(u_1\otimes v_1 + \cdots + u_p\otimes v_p\bigr) + A''\right)V=
 s_1(A)\big( I_p \oplus A_{p+1}\big),
\end{align}
where $I_p = E_{11} + E_{22} + \cdots + E_{pp} = f_1\otimes e_1 + f_2\otimes e_2 + \cdots + f_p\otimes e_p$,
and $A_{p+1}$ is some operator disjoint from $I_p$ with $\|A_{p+1}\|<1$.

When $A$ does not have
a norm attaining vector, that is, $\|A\| > \|Ah\|$ for all unit vector $h\in H$,
there will be sequences $\{v_n\}_n$, $\{u_n\}_n$
of unit vectors in $H,K$, which have no convergent subsequence, such that $\|A\|> \langle Av_n, u_n\rangle \to \|A\|$.
In this case,
$A$ does not have finite peak rank, as $s_1(A)=s_2(A)=\cdots = s_p(A)=\cdots= \|A\|$.
See \cite{P}, especially its \S2.11, for the general theory of $s$-numbers.

It is clear that if $A,B\in B(H,K)$ both have norm attaining vectors such that the form \eqref{eq:UAV} holds for $A$, and
$$
UBV = s_1(B)\big([b]\oplus B_2\big)
$$
with $|b|\geq \|B_2\|$ (resp.\ $b\geq \|B_2\|$),
then $A, B$ are parallel (resp.\ form a TEA pair). However, there are cases in which for a parallel
pair $A,B$ in $B(H,K)$, not both $A$ and $B$ have norm attaining vectors. We have
the following criteria from Zamani and Moslehian \cite{ZM15}; in case when
$H$ or $K$ has finite dimension, all $\psi_n$ can be chosen to be a fixed unit vector $\psi$.

\begin{lemma}[{\cite[Theorem 3.3]{ZM15}}]\label{lem:ZM15}
Two operators $A,B\in B(H,K)$ form a parallel (resp.\ TEA) pair if and only if
there is a sequence $\{\psi_n\}$
of unit vectors in $H$ such that $|\langle A\psi_n,B\psi_n\rangle|\to \|A\|\|B\|$
(resp.\ $\langle A\psi_n,B\psi_n\rangle\to \|A\|\|B\|$).
\end{lemma}

Denote by $F(H,K)$ the vector subspace of $B(H,K)$ consisting of all finite rank operators.
Let us call
a vector subspace $S(H,K)$ of $B(H,K)$ \emph{standard}
if it contains $F(H,K)$; and it is
\emph{unitarily invariant}
if $US(H,K)V= S(H,K)$ for all unitary operators $U$ on $K$ and $V$ on $H$.

\begin{lemma}\label{lem:standardOS}
  Every nonzero unitarily invariant vector subspace $S(H,K)$ of $B(H,K)$ is standard.
\end{lemma}
\begin{proof}
  Let $A\in S(H,K)$ be nonzero. Then there is a unit vector $e\in H$ such that $Ae\neq 0$. Consider the unitary
  operator $V_e=I-2e\otimes e$ on $H$. Since $S(H,K)$ is unitarily invariant,
  the rank one operator $B=A-AV_e= 2Ae\otimes e\in S(H,K)$. Consequently,
  $$
  F(H,K)= \operatorname{span} \{UBV:   U\in B(K,K), V\in B(H,H) \ \text{are unitary operators}\}\subseteq S(H,K).
  $$
\end{proof}

In the following, we assume both $H,K$ are (real or complex) Hilbert spaces of
 dimension at least $2$. Let $S(H,K), S_1(H,K), S_2(H,K)$ denote
  unitarily invariant standard vector subspaces of $B(H,K)$.
Examples of such subspaces include finite rank operators $F(H,K)$, compact operators $\mathcal{K}(H,K)$, and
Schatten $p$-class operators $\mathcal{S}_p(H,K)$ for $1\leq p \leq \infty$.

The following result is an important tool for many preserver problems.

\begin{lemma}[{Hou \cite{Hou}}]\label{lem:Hou}
  Let $T:S_1(H,K)\to S_2(H,K)$ be a linear map sending rank-one
operators to rank-one operators and disjoint rank-one operators to disjoint
rank-one operators. Then either
  \begin{enumerate}
  \item there are injective linear maps $B: K\to K$
  and $C:H\to H$ such that $T(k\otimes h) = Bk\otimes Ch$, or
  \item there are injective conjugate linear maps $B: H\to K$ and $C: K\to H$ such that $T(k\otimes h) = Bh\otimes Ck$,
  \end{enumerate}
  for every rank one operator $k\otimes h$ from $H$ into $K$.

In either case, both $B,C$ are scalar multiples of linear or
conjugate-linear isometric embeddings, respectively.
\end{lemma}
\begin{proof}
  The rank-one preserver classification also allows degenerate maps
whose rank-one images all have a common range line or a common initial
line. Such maps are excluded here: the domain contains two nonzero disjoint
rank-one operators, and their images are nonzero and disjoint. The remaining
two forms follow from Hou \cite[Lemmas 1.1 and 1.2]{Hou};
  although the statements and the proofs are stated for rank one linear preservers of $B(X)$ for a
complex Banach space $X$ in \cite{Hou}, they are
also valid and easily adapted to our case. Note that we do not assume any continuity of $T$, as
we do not assume any continuity on the linear maps $B, C$, either.

For the second part,
if $h, h'$ are orthogonal vectors in $H$ and $k, k'$ are orthogonal vectors
 in $K$ then
 $(k\otimes h)^*(k'\otimes h') = (h\otimes k)(k'\otimes h')= \langle k',k\rangle (h\otimes h')=0$,
 and $(k\otimes h)(k'\otimes h')^*= \langle h',h\rangle (k\otimes k')=0$. In other words,
 $k\otimes h$ and $k'\otimes h'$ are disjoint rank one operators in $F(H,K)$,
and thus so are their images under $T$. In case (1), $T(k\otimes h)=Bk\otimes Ch$ is disjoint
from $T(k'\otimes h')=Bk'\otimes Ch'$; namely,
 both $Bk, Bk'$ and $Ch, Ch'$ are orthogonal pairs. Thus the linear maps $B,C$
 send orthogonal vectors to orthogonal vectors. For any unit vectors $h,h'\in H$, there exists $e^{i\theta}$ such that $e^{-i\theta}\langle h,h'\rangle=|\langle h,h'\rangle|$.
Then $h+e^{i\theta}h'$ and $h-e^{i\theta}h'$ are orthogonal, so are $Ch+e^{i\theta}Ch'$ and $Ch-e^{i\theta}Ch'$.
 This implies $\|Ch\|=\|Ch'\|$. It follows that $C$ is a scalar multiple of an injective isometry; the same argument applies to $B$.
 In a similar manner, we can verify case (2), and complete the proof.
\end{proof}

\begin{lemma}\label{lem:finite-infinite}
  Assume that the dimension of either $H$ or $K$ is infinite and the other is finite.
 Any bijective linear map $T:B(H,K)\to B(H,K)$ preserving
parallel pairs (or TEA pairs) in both directions assumes
the form $A\mapsto \gamma UAV$ for
  a positive scalar $\gamma$, and some unitary operators $U\in B(K)$ and $V\in B(H)$.
\end{lemma}
\begin{proof}
If $T$ preserves TEA pairs in both directions, linearity and
scalar phase adjustment imply that it also preserves parallel
pairs in both directions. We therefore work with parallel pairs.

By passing to the bijection
$A\mapsto (T(A^\TT))^\TT$ on $B(K,H)$
if necessary, we may assume that
$2\leq \dim K=m<\infty$ and $\dim H=\infty$.
In particular, $B(H,K)=F(H,K)$, and every nonzero operator
has finite peak rank.

For every $A\in B(H,K)$, preservation of parallelism in both
directions gives
\[
 T\bigl(\mathcal P(A)\bigr)=\mathcal P(TA),
 \qquad
 T\bigl(\operatorname{span}\mathcal P(A)\bigr)
 =\operatorname{span}\mathcal P(TA).
\]
Consequently, the two spans have the same codimension.
As in Lemma~\ref{lem:Pk1}, if $A$ has peak rank greater than
one, then $\operatorname{span}\mathcal P(A)$ has codimension
at most one. If $A$ has peak rank one, its parallelism span
has infinite codimension. Hence $A$ has peak rank one if
and only if $TA$ does.

We can now apply the rank-one and disjointness arguments
from the finite rectangular matrix case. For $m\geq3$, use
the proof of $(a)\Rightarrow(d)$ in
Theorem~\ref{thm:Sum-R}; for $m=2$, use Case~2 in the proof
of \cite[Theorem~1.1]{LTWW-matrix}.
The required preservation of peak rank one has just been
established, and the equality of parallelism spans used
in those arguments follows directly from the displayed
identity above, without any appeal to equality of finite
codimensions. Moreover, every operator here has a finite
singular value decomposition, so the finite-rank
decompositions in those arguments remain available.
It follows that $T$ sends disjoint rank-one operators to
disjoint rank-one operators.

By Lemma~\ref{lem:Hou}, there are injective linear maps
$B:K\to K$ and $C:H\to H$, each a scalar multiple of an
isometric embedding, such that
\[
 T(k\otimes h)=Bk\otimes Ch
 \qquad(k\in K,\ h\in H).
\]
The alternative form in that lemma is impossible, since
it would require an injective conjugate-linear map from
the infinite-dimensional space $H$ into the
finite-dimensional space $K$.

Write $B=\alpha U$ and $C=\beta W$, where $\alpha,\beta>0$
and $U:K\to K$, $W:H\to H$ are linear isometric embeddings.
Since $K$ is finite-dimensional, $U$ is unitary.
Every $A\in B(H,K)$ has a finite singular value decomposition,
and therefore linearity yields
\[
 T(A)=BAC^*=\alpha\beta\,UAW^*
 \qquad(A\in B(H,K)).
\]

It remains to show that $W$ is surjective. For every $A$,
\[
 \operatorname{ran}T(A)^*\subseteq W(H).
\]
Given any nonzero $h\in H$, choose a unit vector $k\in K$.
By surjectivity of $T$, the rank-one operator $k\otimes h$
belongs to its range. Thus
\[
 \mathbb Fh
 =\operatorname{ran}(k\otimes h)^*
 \subseteq W(H).
\]
Hence $W(H)=H$, and $W$ is unitary.
Setting $\gamma=\alpha\beta>0$ and $V=W^*$ gives
\[
 T(A)=\gamma UAV
 \qquad(A\in B(H,K)),
\]
as required.
\end{proof}

\begin{theorem}\label{main-2} Assume that $S_1(H,K), S_2(H,K)$ are both unitarily invariant
standard subspaces of
$B(H,K)\neq \bM_2(\IF)$.
Let $T: S_1(H,K)\rightarrow S_2(H,K)$ be a
surjective linear map. Assume that either
\begin{itemize}
\item $TF(H,K)\subseteq F(H,K)$, and $T$ is continuous with respect to the weak operator topology (WOT), or
\item $TF(H,K)= F(H,K)$ and $S_1(H,K)\subseteq \mathcal{K}(H,K)$.
\end{itemize}
Then the following conditions are equivalent.
\begin{enumerate}[{\rm (a)}]
\item
$T(A)$ and $T(B)$ are parallel if and only if $A$ and $B$ are parallel,
 $\forall A, B \in S_1(H,K)$.
\item
$T(A)$ and $T(B)$ form a TEA pair if and only if $A$ and $B$ form a TEA pair, $\forall A, B \in S_1(H,K)$.
 \item
$T$ is a scalar multiple of an isometry, and assumes either the form
\begin{quote}
{\rm (i)}\ $A\mapsto \gamma UAV$, \quad or \quad {\rm (ii)}\  $A \mapsto \gamma UA^\TT V$,
\end{quote}
where $\gamma$ is a positive scalar, and $U,V$ are unitary operators.
In case (ii), $H$ and $K$ have the same dimension.
\end{enumerate}
\end{theorem}

We need some lemmas to prove Theorem \ref{main-2}.

Since the matrix case is completely done in Section \ref{s:fc}, together with Lemma \ref{lem:finite-infinite}
\begin{quote}
 we can assume both Hilbert spaces $H, K$ have infinite dimensions.
\end{quote}
With the following Lemma \ref{lem:inj}, we can further assume that $T$ is bijective in Theorem \ref{main-2}.

\begin{lemma}\label{lem:inj}
  Let $T: S_1(H,K)\to S_2(H,K)$ be a linear map such that $A, B$ form a parallel pair (or TEA pairs) whenever $TA, TB$ do. Then $T$ is injective.
\end{lemma}
\begin{proof}
  Let $TA=0$. Since $TA$ forms a parallel/TEA pair with $TB$, we have $A$ forms a parallel/TEA pair with $B$
   for any $B\in S_1(H,K)$.

 Suppose $Ae\neq 0$ for some unit vector $e\in H$.
  Consider the rank one operator $B=f\otimes e$ where $f$ is a unit vector in $K$ orthogonal to $Ae$.
  Since $A,B$ form a parallel/TEA pair, Lemma~\ref{lem:ZM15}
  gives a sequence of unit vectors $\psi_n$ such that
  $|\langle A\psi_n,B\psi_n\rangle|\to\|A\|$.
  The estimate
  $|\langle A\psi_n,B\psi_n\rangle|
  \leq\|A\|\,|\langle\psi_n,e\rangle|$
  implies $|\langle\psi_n,e\rangle|\to1$.
  After multiplying each $\psi_n$ by a suitable unimodular scalar,
  we therefore have $\psi_n\to e$. Consequently,
  \[
  0=|\langle Ae,f\rangle|
   =|\langle Ae,Be\rangle|=\|A\|,
  \]
  a contradiction. Thus $T$ is injective.
\end{proof}

We use the notation
$$
\mathcal{P}(A) = \mathcal{P}_{S(H,K)}(A) =\{X \in S(H,K):  X \hbox{ is parallel with } A\}.
$$

\begin{lemma}\label{lem:rank1-included} \label{lem1}
Let $H,K$ be Hilbert spaces over $\IF=\IR$ or $\IC$ with orthonormal
bases $\{e_j\}_j, \{f_i\}_i$, respectively, and let
$S(H,K)$ be a unitarily invariant standard vector subspace of $B(H,K)$.
For any $A\in S(H,K)$, one of the following three cases happens.
\begin{enumerate}[(a)]
\item $s_1(A)> s_2(A)$.
In this case, there are unitary operators $V\in B(H)$ and $U\in B(K)$ such that
\begin{equation}\label{eq:S1-form1}
UAV = aE_{11} + A_1,
\end{equation}
in which $a=\|A\|$ and $A_1$ is an operator disjoint from $E_{11}=f_1\otimes e_1$
with $\|A_1\| < a$.
Consequently,
\begin{align}\label{eq:not-S1-form}
\operatorname{span} \mathcal{P}(A) &=
\operatorname{span} \{B\in S(H,K): B\hbox{ is parallel with } A\} \notag\\
&= U^*\Bigl\{bE_{11} + B_1:
b\in \mathbb{F}, B_1\in S(H,K)\ \text{disjoint from $E_{11}$}\Bigr\}V^*.
\end{align}

\item $\IF=\IR$ and $s_1(A)=s_2(A)>s_3(A)$.
In this case, there are unitary operators $V\in B(H)$ and $U\in B(K)$ such that
\begin{equation}\label{eq:S1-form2}
UAV = a(E_{11} + E_{22}) + A_2,
\end{equation}
in which $a=\|A\|$ and $A_2$ is an operator disjoint from
$E_{11} + E_{22}=f_1\otimes e_1 + f_2\otimes e_2$
with $\|A_2\| < a$. Consequently,
 $U\bigl(\operatorname{span} \mathcal{P}(A)\bigr)V$
is exactly the subspace of $S(H,K)$ defined by the single equation
$\langle Be_2,f_1\rangle=\langle Be_1,f_2\rangle$. In particular,
$\operatorname{span}\mathcal P(A)$ has codimension one.

\item In all other cases, $\operatorname{span} \mathcal{P}(A)$ is WOT dense in $S(H,K)$.
\end{enumerate}
\end{lemma}
\begin{proof}
The case $A=0$ is immediate. Normalize $\|A\|=1$ otherwise.
If $s_1(A)>s_2(A)$, the isolated maximal singular value gives
\eqref{eq:S1-form1}. An operator parallel with $A$ has the same norming
initial and range lines. Thus it has the indicated block form, with the
modulus of its scalar corner dominating the norm of its other corner.
Taking differences after increasing that scalar corner gives
\eqref{eq:not-S1-form}.

Suppose next that $A$ has at least two orthogonal norming initial vectors.
By unitary changes of coordinates, assume $Ae_i=f_i$, $i=1,2$.
For arbitrary $B\in S(H,K)$ disjoint from $E_{ii}$ and sufficiently large
$t>0$, both $tE_{ii}+B$ and $tE_{ii}$ are parallel with $A$.
Their difference shows that $B\in\operatorname{span}\mathcal P(A)$.
Together with $E_{ii}$, these operators span all of $S(H,K)$ except
possibly the two entries $(1,2)$ and $(2,1)$: rows and columns involving
the first two coordinates have finite rank and belong to $S(H,K)$.
The operator $E_{12}+E_{21}$ is also parallel with $A$.
Over $\IC$, so is $\mathrm i(E_{12}-E_{21})$, giving the entire space.
Over $\IR$, a third orthogonal norming vector likewise gives the entire
space. If there are exactly two norming vectors and a spectral gap below
one, a parallel operator has a common norming initial vector with $A$
in $\operatorname{span}\{e_1,e_2\}$. Its leading real $2\times2$ corner
has this vector as both a right and a left singular vector with the same
signed singular value, hence that corner is symmetric. This proves the
codimension-one description in (b).

It remains to treat the cases in which the preceding arguments do not
apply. If the norm is attained at exactly one or two independent initial
vectors but there is no gap, let $P$ be the projection onto
$\ker(I-A^*A)$ and put $A_0=A(I-P)$. If the norm is not attained, put
$P=0$ and $A_0=A$. Then $P$ has finite rank,
$A_0\in S(H,K)$, $\|A_0\|=1$, and $A_0$ does not attain its norm.
Moreover, $A$ and $A_0$ agree on $(I-P)H$.

Write $A_0=W|A_0|$ for its polar decomposition, and choose $0<b_n\uparrow1$.
Let
\[
 P_n=\mathrm E_{A_0^*A_0}([0,b_n)),\qquad
 Q_n=I-W\mathrm E_{A_0^*A_0}([b_n,1])W^*.
\]
Both are projections and $P_n\to I$, $Q_n\to I$ strongly, because
$\mathrm E_{A_0^*A_0}(\{1\})=0$.
For fixed unit vectors $e\in H,f\in K$, put
$R_n=(Q_nf)\otimes(P_ne)$. The operator $R_n$ maps $P_nH$ into $Q_nK$,
whereas $A_0(I-P_n)$ maps the orthogonal initial space into
$(I-Q_n)K$. Since $\|A_0P_n\|\leq\sqrt{b_n}<1$, for sufficiently small
$\epsilon_n>0$,
\[
 \|A_0-\epsilon_nR_n\|
 =\max\{\|A_0(I-P_n)\|,\|A_0P_n-\epsilon_nR_n\|\}=1.
\]
Choose unit vectors $h_m$ in spectral subspaces of $A_0^*A_0$
corresponding to intervals approaching one. For fixed $n$ and sufficiently
large $m$, $R_nh_m=0$ and $Ah_m=A_0h_m$, while $\|Ah_m\|\to1$.
Lemma~\ref{lem:ZM15} shows that both $A_0$ and
$A_0-\epsilon_nR_n$ are parallel with $A$. Hence
$R_n\in\operatorname{span}\mathcal P(A)$.
As $R_n\to f\otimes e$ in norm, the WOT closure of this span contains
all finite-rank operators and therefore all of $S(H,K)$.
\end{proof}

\begin{lemma} \label{lem2}  \label{lem:B-diag}
Let $k\geq3$ and let $B_1,\dots,B_k\in S(H,K)$ be such that
$s_1(B_j) > s_2(B_j)$ and $B_i \notin \mathcal{P}(B_j)$ for $i,j \in \{1, \dots, k\}$
with $i\ne j$. Suppose
$\sum_{l=1}^k \mu_l B_l \in \mathcal{P}(B_j)$ whenever $|\mu_j| = \displaystyle\max_l |\mu_l|$.
Then there are unitaries $U\in B(K), V\in B(H)$ and $\gamma > 0$ such that
$UB_jV = \gamma (E_{jj} + \hat{Z}_j)$,
where $\hat{Z}_j$ is disjoint from $I_k:=E_{11}+ \cdots+ E_{kk}$ with $\|\hat{Z}_j\|<1$
 for $j = 1,\dots, k$.
\end{lemma}
\begin{proof}
Choose unit vectors $x_j\in H,y_j\in K$ such that
$B_jx_j=\|B_j\|y_j$ and $B_j^*y_j=\|B_j\|x_j$.
Since $s_1(B_j)>s_2(B_j)$, these norming initial and range lines are unique.
Write $B(\mu)=\sum_l\mu_lB_l$. For $\mu_j=1$ and $|\mu_l|<1$ when
$l\ne j$, the hypothesis and the simple maximal singular value of $B_j$
give
\[
 B(\mu)x_j\in\IF y_j,\qquad B(\mu)^*y_j\in\IF x_j.
\]
Varying each coefficient in this open set shows that every $B_l$ satisfies
these same two inclusions. Thus each pair of lines $(\IF x_j,\IF y_j)$
reduces every $B_l$ in the rectangular sense.

For $i\ne j$, write $B_ix_j=d_{ji}y_j$. Necessarily
$|d_{ji}|<\|B_i\|$, since equality would make $B_i$ parallel with $B_j$.
The vectors $x_i,x_j$ are eigenvectors of $B_i^*B_i$ for distinct
eigenvalues, so they are orthogonal. The same argument applied to
$B_iB_i^*$ gives $y_i\perp y_j$. Consequently, in these orthonormal
coordinates,
\[
 B_l=\operatorname{diag}(d_{1l},\ldots,d_{kl})\oplus Z_l,
 \qquad d_{ll}=\|B_l\|>\max\{\|Z_l\|,|d_{jl}|:j\ne l\}.
\]
Let $D=(d_{jl})$. Whenever $|\mu_j|$ is maximal, the hypothesis gives
$|\langle B(\mu)x_j,y_j\rangle|=\|B(\mu)\|$; in particular the $j$th
coordinate of $D\mu$ has maximal modulus. Applying
Lemma~\ref{lem:D-matrix} to $D^{\TT}$ and using $k\geq3$ yields
$D=\gamma I_k$ for some $\gamma>0$. Dividing each $Z_l$ by $\gamma$
gives the required form.
\end{proof}

\begin{lemma}\label{lem:Spk1}\label{prop:T-1}
Let $T: S_1(H,K)\to S_2(H,K)$ be a surjective linear map such that
$TA$ is parallel with $TB$ if and only if $A$ is parallel with $B$.
If it holds either
\begin{enumerate}[{\rm (a)}]
  \item $TF(H,K)=F(H,K)$, or
  \item $TF(H,K)\subseteq F(H,K)$ and $T$ is WOT continuous,
\end{enumerate}
then $T$ sends orthogonal rank one operators to orthogonal rank one operators.
\end{lemma}
\begin{proof}
By Lemma~\ref{lem:inj}, $T$ is bijective. Parallelism preservation in both
directions gives
\[
 T(\operatorname{span}\mathcal P(A))
 =\operatorname{span}\mathcal P(TA).
\]
For a rank-one partial isometry $E$, the left-hand span before applying
$T$ has infinite codimension. Since $TE$ has finite rank,
Lemma~\ref{lem:rank1-included} implies $s_1(TE)>s_2(TE)$: the other
finite-rank possibilities have span of codimension at most one, as shown
in its proof.

Given any finite family of disjoint rank-one partial isometries $E_j$,
enlarge it, if necessary, to a family of at least three by adding disjoint
rank-one partial isometries in the domain. This is possible because both
Hilbert spaces are infinite-dimensional. The operators $B_j=TE_j$
satisfy Lemma~\ref{lem:B-diag}: the required implication follows from
$\sum_l\mu_lE_l\parallel E_j$ whenever $|\mu_j|=\max_l|\mu_l|$.
It follows that
\[
 \operatorname{rank}(TA)\geq\operatorname{rank}(A)
 \qquad(A\in F(H,K)).
\]
Under (a), the same argument applies to $T^{-1}$, so equality holds.

Under (b), fix a rank-one partial isometry $F$ in the range and put
$C=T^{-1}F$. The span $\operatorname{span}\mathcal P(F)$ is a proper
WOT-closed subspace. WOT continuity of $T$ therefore implies that
$\operatorname{span}\mathcal P(C)$ is not WOT dense. Its codimension is
infinite. Lemma~\ref{lem:rank1-included} excludes both (b) and (c), so
$s_1(C)>s_2(C)$. For a finite disjoint family $F_j$, enlarged to at least
three members, apply Lemma~\ref{lem:B-diag} to $T^{-1}F_j$.
If $TA$ has rank $h$ and $A$ has finite rank $k$, the resulting diagonal
corner shows $k\geq h$. Combined with the previous inequality, this gives
$h=k$ also under (b).

In particular, $TE_j$ has rank one. Its remainder in the representation
of Lemma~\ref{lem:B-diag} must vanish. Hence the images of disjoint
rank-one operators are disjoint rank-one operators.
\end{proof}

\begin{proof}[Proof of Theorem \ref{main-2}]
  The implications (c) $\Rightarrow$ (b) $\Rightarrow$ (a) are clear.
Assume (a). Lemmas~\ref{lem:Hou} and~\ref{lem:Spk1} give a
representation on $F(H,K)$ by scalar multiples of isometric embeddings,
possibly with a transpose. In the WOT-continuous case, this representation
extends to $S_1(H,K)$ by WOT density. The extension is WOT continuous,
including in the transpose case, as is seen from its matrix coefficients.
The ranges and initial spaces of all images lie in the corresponding
closed ranges of the implementing embeddings. Since $T$ is onto a space
containing every rank-one operator, these closed ranges must be the full
Hilbert spaces. The embeddings are therefore surjective and give the
unitaries in (c).

Suppose now that $TF(H,K)=F(H,K)$. We show that the action of $T$ assumes either the form
given in (i) or (ii) in (c).
Here $TF(H,K)=F(H,K)$ makes the implementing isometric embeddings
surjective already on the finite-rank part. By Lemmas~\ref{lem:Hou}
and~\ref{lem:Spk1}, after some unitary transformations, a nonzero scalar multiplication, and maybe a transposition,
we can assume that $TA=A$ for all finite rank operators $A\in F(H,K)$.

Let $A\in S_1(H,K)\subseteq \mathcal{K}(H,K)$ be a compact operator of infinite rank.
We claim that $\|TA\|\leq 3\|A\|$.
Write $A = \sum_j \lambda_j P_j$ as a (norm convergent)
countable linear sum of disjoint finite rank partial
isometries $P_j$ with $\lambda_1>\lambda_2 > \cdots > \lambda_j\downarrow 0$.
Let $P_1 = \sum_{j=1}^k u_j\otimes v_j$ be a \emph{finite} sum of disjoint rank one partial isometries.
Let $A_1=A-\lambda_1\sum_{j=1}^{k-1}u_j\otimes v_j$.
Since $A_1\parallel u_k\otimes v_k$, we have
$TA_1\parallel T(u_k\otimes v_k)=u_k\otimes v_k$.
Parallelism with a rank-one partial isometry gives
\[
 TA_1=\lambda u_k\otimes v_k+B,\qquad
 B\perp u_k\otimes v_k,\qquad \|B\|\leq|\lambda|=\|TA_1\|.
\]
Indeed, a norming sequence for the rank-one operator converges, after
scalar phase adjustments, to its norming initial vector; the equality
case of Cauchy--Schwarz then yields this block decomposition.
No compactness of $TA_1$ or spectral gap for $TA_1$ is required.
If $|\lambda|>2\lambda_1$, then
$T(A-\lambda u_k\otimes v_k)\ne0$ by injectivity (the operator $A$ has
infinite rank). In this case $A-\lambda u_k\otimes v_k$ is parallel with $u_k\otimes v_k$, while
$T(A-\lambda u_k\otimes v_k) = \lambda_1 \sum_{j=1}^{k-1} u_j\otimes v_j + B$ is disjoint from $T(u_k\otimes v_k)=u_k\otimes v_k$.
This contradiction shows that
 $|\lambda|\leq 2\lambda_1$, and $\|TA\| \leq \|TA_1\| + \|\lambda_1\sum_{j=1}^{k-1} u_j\otimes v_j\|\leq
 3\lambda_1=3\|A\|$, as claimed.

Observe that
$$
TA = \lambda_1 P_1 + T(\sum_{j\geq 2} \lambda_j P_j) = \cdots = \sum_{j=1}^r \lambda_j P_j + T(\sum_{j\geq r+1} \lambda_j P_j),
$$
with $\|T(\sum_{j\geq r+1} \lambda_j P_j)\|\leq 3 \|\sum_{j\geq r+1} \lambda_j P_j\| = 3 \lambda_{r+1}\downarrow 0$.
Consequently,
$TA=\sum_{j\geq 1}\lambda_j P_j = A$ for every compact operator $A\in S_1(H,K)$.
Since $S_1(H,K)$ consists of compact operators, $T$ assumes the form given in either (i) or (ii) in (c) as asserted.
\end{proof}

The following examples tell us that the surjectivity assumption on $T$ is indispensable.

\begin{example}\label{eg:nonsurj}
(a) Let $H, K$ be infinite dimensional Hilbert spaces.
Fix unitary identifications $H\otimes H\cong H$ and
$K\otimes K\cong K$. These allow us to regard an operator
$A\otimes C:H\otimes H\to K\otimes K$ as an element of $B(H,K)$.
Let $C \in B(H,K)$ be fixed with $s_1(C) > s_2(C)$.
Define $T(A) = A \otimes C$ for all $A\in B(H,K)$.
It is clear that $T:B(H,K)\to B(H,K)$ preserves parallel/TEA pairs
 in \emph{both directions}; namely,
 $$
 A, B \ \text{form a parallel/TEA pair}\quad\Longleftrightarrow\quad
 TA, TB \ \text{form a parallel/TEA pair.}
 $$
Unlike the finite dimensional case,
$T$ need not assume the form $A\mapsto \gamma UAV$ or $A\mapsto \gamma UA^\TT V$, in general.

(b) Let $H=K$ be an infinite dimensional separable Hilbert space
  with an orthonormal basis $\{e_1, e_2, \ldots\}$.
  Let $S=\sum_{n\geq 1} e_{n+1}\otimes e_n$ be the unilateral shift operator.
  Let $T: B(H,H)\to B(H,H)$ be the linear isometry defined by
  $$
  A \mapsto \langle Ae_1, e_1\rangle e_{1}\otimes e_1 + SAS^*;
  $$
  in terms of infinite matrices,
  $$
  A=(a_{ij}) \mapsto   \begin{bmatrix} a_{11} & 0 \\ 0 & A\end{bmatrix}.
  $$
  Then $T$ is WOT continuous, sends finite rank operators to finite rank operators,
   and sends parallel/TEA pairs to parallel/TEA pairs in both directions. However, the range space of $T$ is not standard or unitarily invariant,
  $T$ does not send (orthogonal) rank one operators to (orthogonal) rank one operators, and $T$ does not assume the form
  $A\mapsto \gamma UAV$ or $A\mapsto \gamma UA^\TT V$ for any unitaries $U,V$.
\end{example}

Recall that a type I factor von Neumann algebra is $B(H)=B(H,H)$ for some complex Hilbert space $H$. When $H$ has finite dimension $n$, we have $B(H)=\bM_n(\IC)$.
While bijective linear parallelism preservers of matrix algebras are completely characterized in
Propositions \ref{prop:matrices}, \ref{prop-parallel-C} and \ref{prop-parallel-R}, the following corollary of
Theorem \ref{main-2} covers the infinite dimensional case.

\begin{corollary}\label{cor:typeIfactors}
Let $H,K$ be nonzero Hilbert spaces over the same field
$\mathbb F=\mathbb R$ or $\mathbb C$, with $\dim H\ne2$.
Let $T:B(H)\to B(K)$ be a WOT-continuous surjective linear map
preserving parallel pairs (or TEA pairs) in both directions, and suppose
that $T(F(H))\subseteq F(K)$.
Then $\dim H=\dim K$, and there are unitaries
$U:H\to K$, $V:K\to H$ and a scalar $\gamma>0$ such that
\[
 T(A)=\gamma UAV\quad\text{for all $A\in B(H)$},
 \qquad\text{or}\qquad
 T(A)=\gamma UA^\TT V\quad\text{for all $A\in B(H)$}.
\]
In particular, $T$ is a scalar multiple of an isometry.
\end{corollary}
\begin{proof}
The one-dimensional case is immediate. Otherwise, the argument of
Lemma~\ref{lem:inj} shows that $T$ is injective, and hence bijective.
If either Hilbert space is finite-dimensional, then both are
finite-dimensional and
$(\dim H)^2=(\dim K)^2$. The conclusion then follows from
Proposition~\ref{prop:matrices}, after a unitary identification.
We therefore assume that both spaces are infinite-dimensional.

If the hypothesis concerns TEA pairs, linearity and scalar phase
adjustment imply preservation of parallel pairs in both directions.
The proof of Lemma~\ref{lem:Spk1}(b) applies equally to maps from
$B(H)$ onto $B(K)$: its arguments concerning parallelism sets,
WOT continuity, and finite-rank operators do not require
$H=K$. Thus $T$ sends disjoint rank-one operators to disjoint
rank-one operators.

The corresponding different-space version of
Lemma~\ref{lem:Hou} now gives a scalar $\gamma>0$ and linear
isometric embeddings $U,V:H\to K$ such that either
\[
 T(A)=\gamma UAV^*
 \qquad(A\in F(H)),
\]
or
\[
 T(A)=\gamma UA^\TT V^*
 \qquad(A\in F(H)).
\]
In the second case, the conjugate-linear implementing maps are
composed with the conjugation associated with the fixed
orthonormal basis of $H$ to obtain the linear maps $U,V$.

Both displayed maps are WOT continuous. Moreover, $F(H)$ is
WOT dense in $B(H)$: if $P_E$ denotes the projection onto a
finite-dimensional subspace $E\subseteq H$, then
$P_E\to I$ strongly as $E$ increases. Indeed, for each $x\in H$,
\[
 \|P_EAP_Ex-Ax\|
 \leq \|A\|\,\|P_Ex-x\|+\|(P_E-I)Ax\|\longrightarrow0.
\]
Thus $P_EAP_E\to A$ strongly, and hence in WOT.
Consequently, WOT continuity of $T$ extends the appropriate
representation to every $A\in B(H)$.

For either representation, every $T(A)$ has range contained
in $U(H)$, and every $T(A)^*$ has range contained in $V(H)$.
Since $T$ is surjective, there is an $A_0\in B(H)$ with
$T(A_0)=I_K$. It follows that
\[
 K=\operatorname{ran}I_K\subseteq U(H),
 \qquad
 K=\operatorname{ran}I_K^*\subseteq V(H).
\]
Thus $U$ and $V$ are surjective isometries, hence unitaries.
In particular, $\dim H=\dim K$, and the asserted conclusion
follows.
\end{proof}

To get a counterpart to Corollary~\ref{cor:typeIfactors}, we note that every simple type I $C^*$-algebra is isomorphic to the algebra $\mathcal{K}(H)$ of compact operators on some Hilbert space $H$.

\begin{corollary}\label{cor:K(H)}
Let $H,K$ be nonzero Hilbert spaces over the same field
$\mathbb F=\mathbb R$ or $\mathbb C$, with $\dim H\ne2$.
Let $T:\mathcal K(H)\to\mathcal K(K)$ be a surjective linear map
preserving parallel pairs (or TEA pairs) in both directions and satisfying
$T(F(H))=F(K)$.
Then $\dim H=\dim K$, and there are unitaries
$U:H\to K$, $V:K\to H$ and a scalar $\gamma>0$ such that
\[
 T(A)=\gamma UAV\quad\text{for all $A\in\mathcal K(H)$},
 \qquad\text{or}\qquad
 T(A)=\gamma UA^\TT V\quad\text{for all $A\in\mathcal K(H)$}.
\]
In particular, $T$ is a scalar multiple of an isometry.
\end{corollary}
\begin{proof}
As in the preceding proof, the one-dimensional case is immediate,
and the argument of Lemma~\ref{lem:inj} shows that $T$ is injective
when $\dim H\geq3$. Thus $T$ is bijective.
The finite-dimensional cases follow as before, so assume that both
$H$ and $K$ are infinite-dimensional.
If the hypothesis concerns TEA pairs, linearity and scalar phase
adjustment show that $T$ preserves parallel pairs in both directions.

The proof of Lemma~\ref{lem:Spk1}(a) applies without requiring
the source and target Hilbert spaces to coincide. Hence $T$ and
$T^{-1}$ preserve disjoint rank-one operators.

Let $\{e_i:i\in I\}$ be an orthonormal basis of $H$.
The operators $T(e_i\otimes e_i)$ are nonzero, rank one, and
pairwise disjoint. Their initial spaces therefore provide an
orthonormal family in $K$, giving $\dim H\leq\dim K$.
Applying the same argument to $T^{-1}$ gives the reverse inequality.
Thus $\dim H=\dim K$.

After identifying $H$ and $K$ by a unitary, the conclusion follows
from Theorem~\ref{main-2}.
\end{proof}

\bigskip
\noindent
{\bf Acknowledgment}

 Li is an affiliate member of the Institute for Quantum Computing, University of Waterloo, and an affiliate faculty of the
Department of Mathematics and Statistics, University of Nevada, Reno, USA. His research was partially supported
by the Simons Foundation Grant 851334.
M.-C. Tsai, Y.-S. Wang and N.-C. Wong are supported by Taiwan NSTC grants 115-2115-M-027-001-MY2, 115-2115-M-005-003-MY2 and 114-2115-M-110-002-MY2,
respectively.

The authors used AI tools (ChatGPT and Gemini) to assist with exposition, language editing, and reviewing selected proofs during revision. The authors remain responsible for all mathematical arguments, references, and conclusions.


\begin{thebibliography}{99}

\bibitem{BL}
L.~B. Beasley and T.~J. Laffey,
\emph{Linear operators on matrices: the invariance of rank-$k$ matrices},
Linear Algebra Appl. \textbf{133} (1990), 175--184.

\bibitem{BCM19}
T.~Bottazzi, C.~Conde, M.~S. Moslehian, P.~W\'{o}jcik, and A.~Zamani,
\emph{Orthogonality and parallelism of operators on various Banach spaces},
J. Aust. Math. Soc. \textbf{106} (2019), no.~2, 160--183.


\bibitem{BFGMP08}
M.~Burgos, F.~J. Fern\'{a}ndez-Polo, J.~J. Garc\'{e}s, J.~M. Moreno, and A.~M. Peralta,
\emph{Orthogonality preservers in $C^*$-algebras, $JB^*$-algebras and $JB^*$-triples},
J. Math. Anal. Appl. \textbf{348} (2008), no.~1, 220--233.

\bibitem{Hou}
J.-C. Hou,
\emph{Rank-preserving linear maps on $B(X)$},
Sci. China Ser. A \textbf{32} (1989), no.~8, 929--940.

\bibitem{Kad51}
R.~V. Kadison,
\emph{Isometries of operator algebras},
Ann. of Math. (2) \textbf{54} (1951), 325--338.

\bibitem{Kan}
C.-H. Kan,
\emph{Ergodic properties of Lamperti operators},
Canad. J. Math. \textbf{30} (1978), no.~5, 1206--1214.

\bibitem{Ket}
B.~Kuzma, C.-K. Li, E.~Poon, and S.~Singla,
\emph{Linear maps preserving parallel pairs with respect to the Ky-Fan $k$-norm},
Linear Algebra Appl. \textbf{687} (2024), 68--90.

\bibitem{Li0}
C.-K. Li,
\emph{Matrices with some extremal properties},
Linear Algebra Appl. \textbf{212/213} (1995), 255--267.

\bibitem{LTWW}
C.-K. Li, M.-C. Tsai, Y.-S. Wang, and N.-C. Wong,
\emph{Nonsurjective maps between rectangular matrix spaces preserving disjointness, triple products, or norms},
J. Operator Theory \textbf{83} (2020), no.~1, 27--53.

\bibitem{LTWW-matrix}
C.-K. Li, M.-C. Tsai, Y.-S. Wang, and N.-C. Wong,
\emph{Linear maps on matrices preserving parallel pairs},
Banach J. Math. Anal. \textbf{19} (2025), Art. 51.

\bibitem{LTWW-lp}
C.-K. Li, M.-C. Tsai, Y.-S. Wang, and N.-C. Wong,
\emph{Linear preservers of parallel/TEA vectors in $L_p(\mu)$ spaces},
J. Math. Anal. Appl. \textbf{553} (2026), no.~2, Art. 129889.

\bibitem{Maligranda08}
L.~Maligranda,
\emph{Some remarks on the triangle inequality for norms},
Banach J. Math. Anal. \textbf{2} (2008), no.~2, 31--41.

\bibitem{Manna27}
J.~Manna, K.~Mandal, K.~Paul, and D.~Sain,
\emph{A study on preservation of parallel pairs and triangle equality attaining pairs},
Linear Algebra Appl. \textbf{753} (2027), 1--21.

\bibitem{NT02}
R.~Nakamoto and S.-E. Takahasi,
\emph{Norm equality condition in triangle inequality},
Sci. Math. Jpn. \textbf{55} (2002), no.~3, 463--466.

\bibitem{P}
A.~Pietsch,
\emph{Eigenvalues and $s$-Numbers},
Cambridge Studies in Advanced Mathematics, vol. 13, Cambridge University Press, Cambridge, 1987.

\bibitem{Popovici12}
D.~Popovici,
\emph{Norm equalities in pre-Hilbert $C^*$-modules},
Linear Algebra Appl. \textbf{436} (2012), no.~1, 59--70.

\bibitem{Wojcik}
P.~W\'{o}jcik,
\emph{Norm-parallelism in classical $M$-ideals},
Indag. Math. (N.S.) \textbf{28} (2017), no.~2, 287--293.

\bibitem{ZM15}
A.~Zamani and M.~S. Moslehian,
\emph{Exact and approximate operator parallelism},
Canad. Math. Bull. \textbf{58} (2015), no.~1, 207--224.

\bibitem{ZM16}
A.~Zamani and M.~S. Moslehian,
\emph{Norm-parallelism in the geometry of Hilbert $C^*$-modules},
Indag. Math. (N.S.) \textbf{27} (2016), no.~1, 266--281.

\end{thebibliography}
\end{document}